\title{Finitely forcible graphons}
\date{March 2011}

\documentclass{article}
\usepackage{amssymb,amsmath,hyperref}
\usepackage{theorem,bbm}
\usepackage{srcltx}

\nonstopmode

\long\def\killtext#1{}
\newtheorem{theorem}{Theorem}[section]
\newtheorem{prop}[theorem]{Proposition}
\newtheorem{lemma}[theorem]{Lemma}
\newtheorem{claim}[theorem]{Claim}
\newtheorem{corollary}[theorem]{Corollary}

\theorembodyfont{\rmfamily}
\newtheorem{example}[theorem]{Example}
\newtheorem{conjecture}{Conjecture}
\newtheorem{remark}[theorem]{Remark}
\newtheorem{definition}[theorem]{Definition}
\newtheorem{question}[conjecture]{Question}

\newenvironment{proof}{\medskip\noindent{\bf Proof. }}{\hfill$\square$\medskip}
\newenvironment{proof*}[1]{\medskip\noindent{\bf Proof of #1. }}{\hfill$\square$\medskip}

\begin{document}

\addtolength{\baselineskip}{3pt} \setlength{\oddsidemargin}{0.2in}

\def\R{{\mathbb R}}
\def\Q{{\mathbb Q}}
\def\N{{\mathbb N}}
\def\Z{{\mathbb Z}}
\def\C{{\mathbb C}}
\def\hom{{\rm hom}}
\def\inj{{\rm inj}}
\def\surj{{\rm sur}}
\def\Inj{{\rm Inj}}
\def\neighb{{\rm neighb}}
\def\inj{{\rm inj}}
\def\ind{{\rm ind}}
\def\sur{{\rm sur}}
\def\PAG{{\rm PAG}}
\def\eps{\varepsilon}
\def\Ge{{\mathbf G}}
\def\Ha{{\mathbf H}}
\def\Aa{{\mathbf A}}
\def\CUT{\text{\rm CUT}}
\def\IO{{\infty\to1}}
\def\GR{\text{\rm GR}}
\def\CLQ{\text{\rm CLQ}}
\def\LT{{\text{\rm LEFT}}}
\def\RT{{\text{\rm RIGHT}}}
\def\Haus{{\rm Hf}}

\def\maxcut{{\sf maxcut}}
\def\MaxCut{{\sf MaxCut}}
\def\wt{\widetilde}
\def\irreg{{\rm irreg}}

\def\lunl{[\hskip-1pt[}
\def\runl{]\hskip-1pt]}

\def\tv{{\rm tv}}
\def\tr{{\rm tr}}
\def\cost{\hbox{\rm cost}}
\def\val{\hbox{\rm val}}
\def\rk{{\rm rk}}
\def\diam{{\rm diam}}
\def\Ker{{\rm Ker}}
\def\QG{{\cal QG}}
\def\QGM{{\cal QGM}}
\def\CD{{\cal CD}}
\def\Pr{{\sf P}}
\def\E{{\sf E}}
\def\Var{{\sf Var}}
\def\Ent{{\sf Ent}}
\def\T{{^\top}}
\def\PD{{\sf Pd}}
\def\SYM{{\sf Sym}}
\def\Lab{{\rm Lab}}
\def\one{{\mathbbm1}}

\def\AA{{\cal A}}\def\BB{{\cal B}}\def\CC{{\cal C}}
\def\DD{{\cal D}}\def\EE{{\cal E}}\def\FF{{\cal F}}
\def\GG{{\cal G}}\def\HH{{\cal H}}\def\II{{\cal I}}
\def\JJ{{\cal J}}\def\KK{{\cal K}}\def\LL{{\cal L}}
\def\MM{{\cal M}}\def\NN{{\cal N}}\def\OO{{\cal O}}
\def\PP{{\cal P}}\def\QQ{{\cal Q}}\def\RR{{\cal R}}
\def\SS{{\cal S}}\def\TT{{\cal T}}\def\UU{{\cal U}}
\def\VV{{\cal V}}\def\WW{{\cal W}}\def\XX{{\cal X}}
\def\YY{{\cal Y}}\def\ZZ{{\cal Z}}

\author{L.~Lov\'asz\footnote{Research supported by OTKA Grant No.~67867.}, B.~Szegedy}
\maketitle

\begin{abstract}
We investigate families of graphs and graphons (graph limits) that
are determined by a finite number of prescribed subgraph densities.
Our main focus is the case when the family contains only one element,
i.e., a unique structure is forced by finitely many subgraph
densities. Generalizing results of Tur\'an, Erd\H{o}s--Simonovits and
Chung--Graham--Wilson, we construct numerous finitely forcible
graphons. Our examples exhibit a new phenomenon, namely that graphons
defined on infinite topological spaces come up naturally in extremal
graph theory. We also give some necessary conditions for forcibility,
which imply that finitely forcible graphons are ``rare'', and exhibit
simple and explicit non-forcible graphons.
\end{abstract}

\tableofcontents

\section{Introduction}

We define the {\it density} $t(F,G)$ of a simple graph $F$ in a
simple graph $G$ as the probability that a random map $V(F)\to V(G)$
is a graph homomorphism, i.e., it preserves edges.

A classical theorem by Tur\'an implies that if a simple graph $G$ has
edge density larger than $1-\frac1k$, then it contains a complete
$(k+1)$-graph; furthermore, if the edge density is $1-\frac1k$ and
the density of complete $(k+1)$-graphs is $0$, then $G$ is a complete
$k$-partite graph with equal color classes. Here two subgraph
densities force a unique structure on a graph $G$. Stability theorems
(Erd\H{o}s and Simonovits \cite{Sim,ESim}) imply that if the
densities are ``close'' to the above values, then the structure of
the graph is ``close'' to the complete $k$-partite graph.

Another interesting theorem of this type is by Chung, Graham and
Wilson \cite{CGW} asserting that if the edge density of $G$ is
``close'' to $1/2$ and the 4-cycle density is ``close'' to $1/16$,
then $G$ is quasi-random, which means (among many other nice
properties) that then the density of an arbitrary fixed graph $F$ is
``close'' to $2^{-|E(F)|}$.

The second theorem is different from the first one in two important
ways. First, this pair of subgraph densities can never be attained by
finite graphs (they can be approximated with arbitrary precision).
Second, stability makes less sense in this setting since
quasirandomness is a less rigid structure. Motivated by their
results, Chung, Graham and Wilson introduced the notion of a forcing
family, which is any set of graphs that can be used to force
quasi-randomness in a similar way. They ask which graph families are
forcing families.

Our paper goes in a slightly different direction. Instead of asking
which graph families can be used to force quasi-randomness, we ask
which structures can be forced by prescribing the densities of
finitely many subgraphs. For this reason we will define forcing
families more generally.

Most of the time we consider finite simple graphs, i.e., graphs
without loops and multiple edges; where we allow multiple edges, we
emphasize this by talking about multigraphs (we never need loops).

\begin{definition}\label{forcdef}
Let $F_1,F_2,\dots,F_k$ be simple graphs and $a_1,a_2,\dots,a_k$ be
real numbers in $[0,1]$. We say that the set
$\{(F_i,a_i):~i=1,\dots,k\}$ is a {\bf forcing family} if there is a
sequence of simple graphs $\{G_n\}_{n=1}^\infty$ with
$\lim_{n\to\infty} t(F_i,G_n)=a_i$ for $1\leq i\leq k$, and for every
such graph sequence
$\lim_{n\to\infty} t(F,G_n)$ exists for every simple graph $F$.
\end{definition}

We note that if $G'$ is obtained from $G$ by replacing every node with
the same number of twin nodes, then $t(F,G')=t(F,G)$ for every simple
graph $F$. Hence in the definition above, we could restrict our attention
to graph sequences with $\lim_{n\to\infty} |V(G_n)|=\infty$.

This definition immediately implies that there is a simple graph
parameter (a function $r:~\FF\mapsto [0,1]$ on the set $\FF$ of
finite simple graphs) such that $\lim_{n\to\infty}t(F,G_n)=r(F)$
whenever $\{G_n\}_{n=1}^\infty$ satisfies $\lim_{n\to\infty}
t(F_i,G_n)=a_i$ for $1\leq i\leq k$. The graph parameter $r$ encodes
the unique structure which is forced.

Let $C_n$, $K_n$ and $P_n$ denote the cycle, complete graph and path
with $n$ nodes, respectively. The result of Chung, Graham and Wilson
mentioned above says in this language that $\{(K_2,1/2),(C_4,1/16)\}$
is a forcing family. The graph parameter describing the limit is
$r(F)=2^{-|E(F)|}$. Simple graph sequences satisfying the conditions
in the definition are called {\it quasirandom}.

The forced structure in definition \ref{forcdef} is best described as
the limit of a graph sequence, using the newly developed theory of
convergent graph sequences.

Let $\WW$ denote the set of bounded symmetric measurable functions of
the form $W:~[0,1]^2\rightarrow \R$, and let $\WW_0\subset\WW$
consist of those functions with range in $[0,1]$. The elements of
$\WW$ are called {\it graphons}.

A graphon $W$ is a {\it stepfunction} if there is a partition
$\{S_1,\dots,S_n\}$ of $[0,1]$ into measurable sets such that $W$ is
constant on each product set $S_i\times S_j$. Every graph $G$ can be
represented by a function $W_G\in\WW_0$: Let $V(G)=\{1,\dots,n\}$.
Split the interval $[0,1]$ into $n$ equal intervals $J_1,\dots,J_n$,
and for $x\in J_i, y\in J_j$ define $W_G(x,y)=\one_{ij\in E(G)}$.

The notion of subgraph densities can be extended to graphons: for a
graph $F=(V,E)$ and graphon $W\in\WW$, we define
\begin{equation} \label{tFW-def}
t(F,W)=\int\limits_{[0,1]^V} \prod_{ij\in E} W(x_i,x_j)\,\prod_{i\in
V}dx_i\,.
\end{equation}
(This definition is meaningful for multigraphs $F$ too.) These
quantities were introduced in \cite{LSz1} and it was proved that if
in a graph sequence $(G_1,G_2,\dots)$ the densities of every fixed
graph form a convergent sequence, then there is a graphon $W\in\WW_0$
such that $t(F,G_n)\to t(F,W)$ for every graph $F$.

If two graphons have the same simple subgraph densities, then they
are called {\it weakly isomorphic} (see Section \ref{GRAPHONS} for
more on this relation). It follows then that also the densities of
multigraphs in them are also equal.

In this paper we reformulate the problem of forcing in terms of
measurable functions. An immediate advantage of this can already be
seen from the simpler definition of finite forcing:

\begin{definition}\label{FINF-A}
Let $\AA\subseteq\WW$. Let $F_1,F_2,\dots,F_k$ be simple graphs and
$a_1,a_2,\dots,a_k\in[0,1]$. We say that the set
$\{(F_i,a_i):~i=1\dots k\}$ is a {\bf forcing family} in $\AA$ if
there is a unique (up to weak isomorphism) graphon $W\in\AA$ with
$t(F_i,W)=a_i$ for every $1\leq i\leq k$. In this case say that $W$
is {\bf finitely forcible} (in $\AA)$, and the family
$\{F_i:~i=1\dots k\}$ is a {\it forcing family for $W$} (in $\AA)$.
\end{definition}

Due to the identity
\begin{equation}\label{EQ:TPROD}
t(F_1F_2,W)=t(F_1,W)t(F_2,W)
\end{equation}
(where $F_1F_2$ denotes the disjoint union of $F_1$ and $F_2$), every
finitely forcible graphon can be forced by a family of connected
graphs.

The two main choices for $\AA$ will be $\AA=\WW$ and $\AA=\WW_0$.
Definition \ref{forcdef} is equivalent with Definition \ref{FINF-A}
when $\AA=\WW_0$.

If a graphon $W\in\WW_0$ is finitely forcible in $\WW$, then it is
also finitely forcible in $\WW_0$, but the reverse is open. A forcing
family for $W$ in $\WW_0$ is not necessarily forcing $W$ in $\WW$.
For example, one can show that the constraints $t(P_3,W)=1/4$ and
$t(C_4,W)=1/16$ force the function $W\equiv 1/2$ among functions in
$\WW_0$, but allows also the function $W\equiv -1/2$ in $\WW$.
(Nevertheless, the graphon $1/2$ is finitely forcible in $\WW$, for
example, by adding the constraint $t(K_2,W)=1/2$.)

Besides the advantage of a simpler definition, the new language
enables us to specify the structure which is forced and to use
analytic methods together with algebraic ones. In this language the
Chung-Graham-Wilson theorem says that if $t(K_2,W)=1/2$ and
$t(C_4,W)=1/16$ for $W\in\WW_0$, then $W$ is the constant $1/2$
function. A generalization of this was proved by Lov\'asz and S\'os
in \cite{LSos}: every stepfunction is finitely forcible in $\WW_0$.

After the Lov\'asz--S\'os result it remained open for a few years
whether only stepfunctions are finitely forcible. This is true in a
one-variable analogue of the forcing problem: a bounded function
$f:~[0,1]\to\R$ is forced by finitely many moments if and only if it
is a stepfunction. In this paper we show that in the 2-variable case
more complicated structures can be forced. One family of these
structures is the indicator function of a level set of a monotone
symmetric 2-variable polynomial. Our other main example has an
iterated (fractal like) structure.

So being a stepfunction is not a characterization of finitely
forcible graphons. The examples mentioned above are, however,
stepfunctions in a weaker sense, i.e., their range is finite. Even
this is not necessary: In Section \ref{OPER} we develop a class of
operations on graphons that preserve finite forcibility in $\WW$, and
applying these to the first type of examples, we construct finitely
forcible graphons whose range is a continuum. In the finite language
this implies the surprising fact that we can create an extremal
problem involving the densities of finitely many subgraphs such that
in the unique asymptotically optimal solution a continuous spectrum
of probabilities for quasi-randomness appears.

One implication of these results could be that they provide
prototypes of possible extremal graphs other than the usual ones
modeled by stepfunctions.

Of course, instead of forcing a given graphon by a finite number of
density constraints, we can also try to force various properties. Let
$\BB\subseteq\AA\subseteq\WW$ be closed under weak isomorphism. We
say that $\BB$ is {\it finitely forcible in $\AA$} if there exists a
family $\{(F_i,a_i):~i=1\dots k\}$ (where the $F_i$ are simple
graphs) and $a_i\in[0,1]$) such that a graphon $W\in\AA$ satisfies
the constraints $t(F_i,W)=a_i$ ($i=1\dots k$) if and only if it is in
$\BB$. While the study of this generalization is not the goal of this
paper, we do need certain facts about forcing some simple properties,
which are discussed in Section \ref{PROPS}.

Since finitely forcible graphons can be described by finitely many
real numbers, we believe that they are very special. However, a full
characterization seems to be very difficult. In fact, it is not easy
to find necessary conditions for finite forcibility. We present a few
in this paper and formulate many open problems in this direction.

\section{Preliminaries}

\subsection{Graphons}\label{GRAPHONS}

Recall that a graphon is a symmetric bounded measurable function
$W:~[0,1]^2\to\R$. Sometimes it is more convenient to use
two-variable functions on other probability spaces than $[0,1]$; this
does not add real generality, however, since a 2-variable function on
any probability space can be replaced by an ``equivalent'' graphon on
$[0,1]$; see \cite{BCL} for more details.

Also recall that two graphons are weakly isomorphic if have the same
simple subgraph densities. We denote by $[W]$ the set of all graphons
weakly isomorphic to $W$. In \cite{BCL} various characterizations of
weakly isomorphic graphons were given, of which we need the
following. For a graphon $W$ and map $\varphi:~[0,1]\to[0,1]$, define
$W^\varphi(x,y)=W(\varphi(x),\varphi(y))$.

\begin{theorem}\label{THM:WEAKISO}
Two graphons $U$ and $W$ are weakly isomorphic if and only if there
are measure preserving maps $\varphi,\psi:~[0,1]\to[0,1]$ such that
$U^\varphi=V^\psi$.
\end{theorem}

Every graphon defines a kernel operator $T_W:~L_1[0,1]\to
L_\infty[0,1]$, by
\[
(T_Wf)(x)=\int\limits_0^1 W(x,y)f(y)\,dy.
\]
We can also consider $T_W$ as an operator $L_\infty[0,1]\to L_1[0,1]$
or $L_2[0,1]\to L_2[0,1]$. in the latter case it is a Hilbert-Schmidt
operator, and hence it has a discrete spectrum
$\{\lambda_1,\lambda_2\dots\}$ such that the eigenvalues tend to $0$
(in particular, every nonzero eigenvalue has finite multiplicity).
Furthermore, it has a spectral decomposition
\[
W(x,y) \sim \sum_k \lambda_k f_k(x)f_k(y),
\]
where $f_k$ is the eigenfunction belonging to the eigenvalue
$\lambda_k\not=0$ with $\|f\|_2=1$. It can be shown that the eigenfunctions
$f_k$ are bounded. The $\sim$ sign indicates that
the series on the right may not be
convergent pointwise, only in $L_2$; but one has
\[
\sum_{k=1}^\infty \lambda_k^2 = \int\limits_{[0,1]^2}
W(x,y)^2\,dx\,dy \le \|W\|^2_\infty.
\]
The spectral decomposition can be used, among others, to compute the
inner product of $W$ with any function $U\in L_2[0,1]^2$:
\[
\int\limits_{[0,1]^2} W(x,y)U(x,y)\,dx\,dy = \sum_k \lambda_k
\int\limits_{[0,1]^2} f_k(x)f_k(y)U(x,y)\,dx\,dy.
\]

The {\it cut-norm} introduced in \cite{FK} is defined for $W\in\WW$
by
\[
\|W\|_\square=\sup_{S,T\subset [0,1]} \Bigl|\int\limits_{S\times T}
W(x,y)\,dx\,dy\Bigr|,
\]
where the supremum goes over measurable subsets of $[0,1]$. This is
within a factor of $4$ to the operator norm
$\|T_W\|_{L_\infty[0,1]\to L_1[0,1]}$.

Let $W_1,W_2,\dots\in\WW$ be graphons (a finite or countable
sequence), and let $a_1,a_2,\dots$ be positive real numbers such that
$\sum_i a_i=1$. We use $\sum_i a_iW_i$ to denote the pointwise linear
combination. We also define the weighted direct sum $W=\bigoplus_i
(a_i)W_i$ as the graphon on $[0,1]$ as follows: we break the $[0,1]$
interval into intervals $J_1,J_2,\dots$ of length $a_1,a_2,\dots$,
take homothetical maps $\phi_i:~J_i\to[0,1]$ and define
$W(x,y)=W_i(\phi_i(x),\phi_i(y))$ if $(x,z)\in J_i\times J_i$ for
some $i$, and $W(x,y)=0$ otherwise. See \cite{Jan} for more on
weighted direct sums and decomposing graphons into connected
components.

Somewhat confusingly, we can introduce three ``product'' operations
on graphons, and we will need all three of them. Let $U,W\in\WW$. We
denote by $UW$ their product as functions, i.e.,
\[
UW(x,y)=U(x,y)W(x,y).
\]
We denote by $U\circ W$ the product of $U$ and $W$ as kernel
operators, i.e.,
\[
(U\circ W)(x,y)=\int\limits_0^1 U(x,z)W(z,y)\,dz.
\]
Finally, we denote by $U\otimes W$ their tensor product; this is
defined as a function $[0,1]^2\times [0,1]^2\to[0,1]$ by
\[
(U\otimes W)(x_1,x_2,y_1,y_2)=U(x_1,y_1)W(x_2,y_2).
\]
However, we can consider any measure preserving map
$\phi:~[0,1]\to[0,1]^2$, and define the graphon
\[
(U\otimes W)^\phi(x,y)=(U\otimes W)(\phi(x),\phi(y)).
\]
These graphons are weakly isomorphic for all $\phi$, and so we can
call any of them the tensor product of $U$ and $W$. We note that the
tensor product has the nice property that $t(F,U\otimes W) = t(F,U)
t(F,W)$ for every graph $F$.

We denote the $n$-th power of a graphon according to these three
multiplications by $W^n$ (pointwise power), $W^{\circ n}$ (operator
power), and $W^{\otimes n}$ (tensor power).

\subsection{Signed, labeled, colored, and quantum graphs}\label{BASIC}

Suppose that the edges of a graph $F$ are partitioned into two sets
$E_+$ and $E_-$. The triple $F'=(V,E_+,E_-)$ will be called a {\it
signed graph}. For $W\in\WW$, we define
\begin{equation} \label{tindFW-def}
t(F',W)=\int\limits_{[0,1]^V} \prod_{ij\in E_+}
W(x_i,x_j)\prod_{ij\in E_-} (1-W(x_i,x_j))\,\prod_{i\in V}dx_i\,.
\end{equation}
If all edges are signed ``$+$'', then $t(F',W)=t(F,W)$. If all edges
are signed ``$-$'', then $t(F',W)=t(F,1-W)$. In general, $t(F',W)$
can be expressed as
\begin{equation}\label{SIEVE}
t(F',W) = \sum_{Y\subseteq E_-} (-1)^{|Y|} t((V,E_+\cup Y), W).
\end{equation}

For a simple graph $F$, let $\widehat F$ denote the complete graph on
$V(F)$ in which the edges of $F$ are signed ``$+$'' and the other
edges are signed ``$-$''. Let $G$ be a simple graph and let $W_G$
denote the associated graphon. Then $t(\widehat{F},W_G)$ is the
probability that a random map $V(F)\to V(G)$ preserves both adjacency
and nonadjacency. If $F$ is fixed and $|V(G)|\to\infty$, then this is
asymptotically $t_\ind(F,G)$, the probability that a random injection
$V(F)\to V(G)$ is an embedding as an induced subgraph.

Let $F=(V,E)$ be a {\it $k$-labeled multigraph}, i.e., a multigraph
with $k$ specified nodes labeled $1,\dots,k$ and any number of
unlabeled nodes. Let $V_0=V\setminus[k]$ be the set of unlabeled
nodes. For $W\in\WW$, we define a function $t_k(F,W):~[0,1]^k\to\R$
by
\[
t_k(F,W)(x_1,\dots,x_k)= \int\limits_{x\in[0,1]^{V_0}}\prod_{ij\in E}
W(x_i,x_j)\,\prod_{i\in V_0} dx_i\,.
\]
Note that $t_0(F,W)=t(F,W)$. We can extend this notation to signed
graphs $\widehat{F}$ to get $t_k(\widehat{F},W)$ in the obvious way.
If all nodes of a signed graph $\widehat{F}=(V,E_+,E_-)$ are labeled,
then
\[
t_k(\widehat{F},W)(x_1,\dots,x_k) =\prod_{ij\in E_+}
W(x_i,x_j)\prod_{ij\in E_-} (1-W(x_i,x_j)).
\]
If a $(k-1)$-labeled graph $F'$ arises from a $k$-labeled graph $F$
by unlabeling node $k$ (say), then
\begin{equation}\label{EQ:UNLABEL}
t_{k-1}(F',W)(x_1,\dots,x_{k-1})=\int\limits_{[0,1]}
t_k(F,W)(x_1,\dots,x_k)\,dx_k.
\end{equation}

A simple but very useful inequality from \cite{LSz1} relates the cut
norm to subgraph densities: for every simple graph $F$ and
$U,W\in\WW$,
\begin{equation}\label{EQ:COUNT}
|t(F,U)-t(F,W)|\le 4
|E(F)|\max\{\|U\|_\infty,\|W\|_\infty\}^{|E(F)|-1} \|U-W\|_\square .
\end{equation}
(One can extend this to simple signed graphs $F$ at the cost of
including $\|1-U\|_\infty$ and $\|1-W\|_\infty$ in the $\max$ on the
right hand side.) We will also use a related (in fact, simpler)
inequality: For all $U,W\in\WW$ and every simple graph $F$ with $k$
nodes we have trivially
\[
\|t_k(F,U)-t_k(F,W)\|_\infty\le
|E(F)|\max\{\|U\|_\infty,\|W\|_\infty\}^{|E(F)|-1} \|U-W\|_\infty,
\]
and hence for all $1\le j\le k$
\begin{equation}\label{EQ:COUNT2}
\|t_j(F,U)-t_j(F,W)\|_\infty \le
|E(F)|\max\{\|U\|_\infty,\|W\|_\infty\}^{|E(F)|-1} \|U-W\|_\infty.
\end{equation}

Let $F_1$ and $F_2$ be two $k$-labeled multigraphs. Their {\it
product} is the $k$-labeled multigraph $F_1F_2$ defined as follows:
we take their disjoint union, and then identify nodes with the same
label (retaining the labels). Clearly this multiplication is
associative and commutative. Note that for $k\geq 2$ the graph
$F_1F_2$ can have multiple edges even if the graphs $F_1$ and $F_2$
are simple. For two $0$-labeled graphs, $F_1F_2$ is their disjoint
union.

For two $k$-labeled graphs $F_1$ and $F_2$ the following
generalization of \eqref{EQ:TPROD} holds:
\begin{equation}\label{EQ:MULTIPL}
t_k(F_1F_2,W)=t_k(F_1,W)t_k(F_2,W)
\end{equation}
where the multiplication on the right hand side is just the pointwise
product of two real functions with the same domain.

A {\it $k$-labeled quantum graph} is a formal finite linear
combination with real coefficients of $k$-labeled multigraphs.
Multigraphs that occur in a $k$-labeled quantum graph with nonzero
coefficient will be called its {\it constituents}. A $0$-labeled
(unlabeled) quantum graph will be called simply a {\it quantum
graph}. Following Razborov \cite{Razb1}, we denote by $\lunl f\runl $
the unlabeled [quantum] graph obtained from the [quantum] graph $f$
by unlabeling every node. A $k$-labeled quantum graph is {\it
simple}, if it is a linear combination of simple graphs. A
$k$-labeled graph is {\it connected}, if every connected component
contains a labeled node. A $k$-labeled quantum graph is {\it
connected}, if every constituent is connected.

We can extend the definition of the product to the product of two
$k$-labeled quantum graphs by distributivity. This way $k$-labeled
quantum graphs form a commutative algebra $\GG_k$; the graph with $k$
nodes, all labeled, and no edge is the unit element in this algebra.

We can define $t_k(f,W):~[0,1]^k\to\R$ for every $k$-labeled quantum
graph $f$ so that it is linear in $f$. Then \eqref{EQ:MULTIPL} will
remain valid.

We can identify a signed graph $\widehat{F}=(V,E_+,E_-)$ with the
quantum graph
\[
\sum_{Y\subseteq E_-} (-1)^{|Y|} (V,E_+\cup Y),
\]
then the two possible definitions of $t(\widehat{F},W)$ give the same
result by \eqref{SIEVE}.

Recall that a graphon $W$ is {\it finitely forcible} (in
$\AA\subseteq\WW$), if $W\in\AA$ and there are a finite number of
simple graphs $F_1,\dots,F_k$ so that whenever a graphon $U\in\AA$
satisfies
\begin{equation}\label{FORCE}
t(F_i,U)=t(F_i,W) \qquad (i=1,\dots,k)
\end{equation}
then $W$ and $U$ are weakly isomorphic. We could be more general and
allow quantum graphs in the forcing family. This would not lead to
more finitely forcible graphons, but we could use forcing constraints
of the special form $t(f_i,W)=0$, where $f_i$ is a quantum graph.

\subsection{Moments of one-variable functions}

We will need a certain one-variable version of finite forcing. Let
$\mathbf{w}=(w_1,\dots,w_r)$, where $w_1,\dots,w_r:~[0,1]\to\R$ are
bounded measurable functions. For $\mathbf{a}=(a_1,\dots,a_r)\in
\N^r$, (where $\N=\{0,1,2,\dots\}$) we define the joint moment of
$\mathbf{w}$ by
\[
M(\mathbf{w},\mathbf{a})=\int\limits_{[0,1]} w_1(x)^{a_1}\cdots
w_r(x)^{a_r}\,dx.
\]
The following theorem follows from classical results of Karlin
\cite{Kar} (see also \cite{Dub}). For any map $\phi:~[0,1]\to[0,1]$
and any $w:~[0,1]\to\R$, we set $w^\phi(x)=w(\phi(x))$.

\begin{theorem}\label{ONE-MOM}
Let $\mathbf{a}_1,\dots,\mathbf{a}_m\in\N^r$, and suppose that
$M(\mathbf{w},\mathbf{a}_j)$ exists for all $j=1,\dots,m$. Then there
is a vector $\mathbf{u}=(u_1,\dots,u_r)$ of stepfunctions such that
\[
M(\mathbf{u},\mathbf{a}_j)=M(\mathbf{w},\mathbf{a}_j)\qquad
(j=1,\dots,m).
\]
\end{theorem}

A certain converse of this theorem is also true:

\begin{prop}\label{ONE-STEP-MOM}
Let $u_1,\dots,u_r:~[0,1]\to\R$ be stepfunctions. Then there is a
finite set of vectors $\mathbf{a}_1,\dots,\mathbf{a}_m\in\N^r$ such
that the values $M(\mathbf{u},\mathbf{a}_j)$ ($j=1,\dots,m$) uniquely
determine the functions $u_i$ up to a measure preserving
transformation of $[0,1]$.
\end{prop}

\begin{proof}
First we prove this for $r=1$:

\begin{claim}\label{CLAIM:STEPMOMENT}
Every stepfunction $u:~[0,1]\to\R$ with $k$ steps is uniquely
determined, up to a measure preserving transformation in the
variable, by its first $2k$ moments.
\end{claim}

Let $u$ have $k$ steps, of sizes $\alpha_1,\dots,\alpha_k$, on which
the value of $u$ is $\beta_1,\dots,\beta_k$, respectively. Then
\[
M(u,a) = \sum_{i=1}^k \alpha_i\beta_i^a.
\]
Let $f:~[0,1]\to\R$ be a function such that $M(u,a)=M(f,a)$ for
$a=1,\dots,2k$. Consider the polynomial
\[
p(x)=\prod_{i=1}^k (x-\beta_i)^2 = \sum_{j=0}^{2k} c_jx^j.
\]
Then
\[
\int\limits_0^1 p(f(x))\,dx = \sum_{j=0}^{2k} c_j M(f,j) =
\sum_{j=0}^{2k} c_j M(u,j) = \int\limits_0^1 p(u(x))\,dx = 0.
\]
Since $p$ is a square, this implies that $p(f(x))=0$ almost
everywhere. Hence $f(x)\in\{\beta_1,\dots,\beta_k\}$ almost
everywhere, i.e., up to a set of measure $0$, $f$ is a stepfunction
attaining the same values as $u$.

Let $\alpha_i'=\lambda(f^{-1}(\beta_i))$ (where
$\lambda$ denotes the Lebesgue measure). Then
\[
\sum_{i=1}^k \alpha_i\beta_i^a=M(u,a) = M(f,a)  = \sum_{i=1}^k
\alpha_i'\beta_i^a \qquad (0\le a\le k-1),
\]
or
\[
\sum_{i=1}^k (\alpha_i-\alpha_i')\beta_i^a = 0\qquad (0\le a\le k-1).
\]
Considering this as a system of linear equations on the differences
$\alpha_i-\alpha_i'$, the determinant of the system is nonzero, which
implies that $\alpha_i-\alpha_i'=0$ for all $i$. So $f$ differs from
$u$ in a measure preserving transformation only.

Now we turn to the general case. Let $u_i$ have $k_i$ steps, and set
$k=k_1k_2\dots k_r$. We claim that if $v_1,\dots,v_r$ is a system of
functions for which
\[
M(\mathbf{v},\mathbf{a}) =M(\mathbf{u},\mathbf{a})
\]
for all $a\in\N^r$ with $\sum_i a_i\le 2k$, then there is a measure
preserving transformation of $[0,1]$ which transforms $v_i$ into
$u_i$ for all $i$.

First, by specifying the first $2k_i$ moments of each $u_i$
separately, we get that $v_i$ is a stepfunction obtained from $u_i$
by a measure preserving transformation of the variable. We want to
argue that we can use the same measure preserving transformation for
every $i$.

Let $\eps>0$ be a sufficiently small real number (smaller than the
minimum difference between two distinct values of any $u_i$, divided
by the maximum of $4r\|u_i\|_\infty$, $i=1,\dots,r$), and consider the function
\[
u=u_1+\eps u_2+\dots +\eps^{r-1} u_r.
\]
Then $u$ is a stepfunction with at most $k$ steps, and the $j$-th
moment of $u$ can be expressed as a polynomial of the joint moments
$M(\mathbf{u},\mathbf{a})$ of the $u_i$ with $\sum_ia_i=j$. Hence the
first $2k$ moments of the function
\[
v=v_1+\eps v_2+\dots +\eps^{r-1} v_r
\]
match the first $2k$ moments of $u$, and so there is a measure
preserving transformation $\phi$ with $v^\phi=u$. By the choice of
$\eps$, the steps of $v_1$ must be mapped by $\phi$ onto the steps of
$u_1$, and so $v_1^\phi=u_1$. But then $(v-v_1)^\phi=u-u_1$, and we
get similarly that $v_2^\phi=u_2$, and so on.
\end{proof}

\subsection{Typical points of graphons}\label{TYPICAL}

We need some technical results about 2-variable functions. In this
section $W$ denotes a measurable function $[0,1]^2\to [0,1]$ (not
necessarily symmetric). Let $R(W)$ denote the set of 1-variable
functions $\{W(x,.):~x\in[0,1]\}$. Clearly $R(W)$ inherits a topology
from $L_1[0,1]$, and it also inherits a probability measure $\pi$
from $[0,1]$.

\begin{definition}
Let $T(W)$ be the set of functions $f\in L_1[0,1]$ such that every
neighborhood of $f$ intersects $R(W)$ in a set with positive measure.
A point $x\in [0,1]$ will be called {\it typical} if $W(x,.)\in T(W)$
and {\it atypical} otherwise.
\end{definition}

\begin{lemma}\label{lem:mes}
Let $W:~[0,1]^2\to [0,1]$ be a measurable function. Then almost every
point of $[0,1]$ is typical.
\end{lemma}

\begin{proof}
If $g \notin T(W)$, then there is an open neighborhood $U_g$ of $g$
in $L_1[0,1]$ such that $\pi(U_g)=0$. Let $U = \bigcup_{g\notin T(W)}
U_g$. Since $L_1[0,1]$ is separable, it contains a countable dense
set $D$, and then every set $U_g$ is the union of all balls with
rational radius centered at a point in $D$ contained in $U_g$. Hence
$U$ equals the union of a countable number of such balls contained in
a $U_g$, $g\notin T(W)$, and hence, it is also the union of a
countable subfamily $\{U_{g_i}:~i\in\N\}$ ($g_i\notin T(W)$). Hence
$\pi(U)=0$. If $x$ is atypical, then $W(x,.) \in U$, and $|\{x:
W(x,.) \in U\}| =\pi(U)=0$.
\end{proof}

\killtext{
\begin{remark}\label{REM:T-GEN}
The proof above can be modified so that instead of $[0,1]$, it works
for graphons defined on any probability space $\Omega$. In this
general case $L_1[\Omega]$ is not necessarily separable, but we can
replace it by the linear space generated by functions $W(x,.)$, which
is separable.
\end{remark}
}

We need the following property of typical points.

\begin{lemma}\label{derand2}
Let $W$ be a graphon, and let $f$ be a $k$-labeled quantum graph such
that the labeled nodes are independent in each multigraph
constituting $f$. Assume that $t_k(f,W)=0$ almost everywhere. Then
$t_k(f,W)(x_1,\dots,x_k)=0$ for every $k$-tuple of typical points.
\end{lemma}

\begin{proof}
We may assume that $\|W\|_\infty \le 1$. Suppose that
$t_k(f,W)(x_1,\dots,x_k)=\eps>0$ with all $x_i$ typical. Let $f=\sum_i \alpha_iF_i$,
$c_f=\sum_i |\alpha_i|\cdot|E(F_i)|$ and $\delta=\eps/(2c_f)$. By the
definition of typical points, there are sets $Z_i\subseteq [0,1]$
with positive measure such that $\|W(x_i,.)-W(z,.)\|_1\le
\delta$ for all $z\in Z_i$. We claim that for every choice of points
$z_i\in Z_i$, we have
\begin{equation}\label{EQ:TKFW}
|t_k(f,W)(x_1,\dots,x_k)-t_k(f,W)(z_1,\dots,z_k)|\le \frac\eps2
\end{equation}
Clearly it suffices to verify this for the case when $f=F$ is a
multigraph. Let $u_1v_1,\dots,u_qv_q$ be the edges of $F$ incident
with the labeled nodes; say, $v_r$ is labeled but $u_r$ is not (here
we use the assumption about $f$). Let $u_{q+1}v_{q+1},\dots,u_mv_m$
be the other edges of $F$, and $U=V\setminus \{1,\dots,k\}$. Using
variables $y_u$ for the unlabeled nodes, we  have
\begin{align*}
|t_k(f,W)(x_1,\dots,x_k)&-t_k(f,W)(z_1,\dots,z_k)|\\
&=\left|\int\limits_{[0,1]^U} \sum_{j=1}^q A_j(x,y)
(W(y_{u_j},x_{v_j})-W(y_{u_j},z_{v_j})) B_j(z,y) \,dy\right|,
\end{align*}
where
\[
A_j(x,y)= \prod_{i<j} W(y_{u_i},x_{v_i})
\]
and
\[
B_j(z,y) = \prod_{j<i\le q} W(y_{u_i},z_{v_i}) \prod_{i>q}
W(y_{u_i},y_{v_i}).
\]
Hence
\begin{align*}
|t_k(f,W)(x_1,\dots,x_k)&-t_k(f,W)(z_1,\dots,z_k)|\\
&\le\sum_{j=1}^q  \int\limits_{[0,1]^U} |A_j(x,y)|\cdot
\bigl|W(y_{u_j},x_{v_j})-W(y_{u_j},z_{v_j})\bigr|
\cdot |B_j(z,y)|\, dy\\
&\le \sum_{j=1}^q \int\limits_{[0,1]^U}
|W(y_{u_j},x_{v_j})-W(y_{u_j},z_{v_j})|\,dy \\
&= \sum_{j=1}^q  \|W(x_{v_j},.)-W(z_{v_j},.)\|_1 \le
q\delta\le\frac{\eps}{2}.
\end{align*}
This proves \eqref{EQ:TKFW}, which in turn implies that
$t_k(f,W)(z_1,\dots,z_k)\not=0$. Since this holds for all $z_i\in
Z_i$, we get that $t_k(f,W)=0$ cannot hold almost everywhere, a
contradiction.
\end{proof}

\begin{remark}\label{REM:TYPICAL}
We can use Lemma \ref{lem:mes} to define a ``normalization'' of
graphons: by modifying a graphon on a set of measure $0$, we can
obtain one in which every point is typical. Lemma \ref{derand2}
implies then that if $t_k(f,W)=0$ almost everywhere, then it is
identically $0$.
\end{remark}

\section{Operations on graphs and graphons}\label{OPER}

We discuss various operations on graphs and graphons in connection
with forcing.

\subsection{Labeling, unlabeling and contraction}

We start with a trivial consequence of \eqref{EQ:UNLABEL}.

\begin{lemma}\label{LEM:UNLABEL0}
Assume that $t_k(g,W)=0$ holds almost everywhere for some $k$-labeled
quantum graph $g$, and let $g_1$ be obtained by unlabeling
$r+1,\dots,k$ in $g$. Then $t_r(g_1,W)=0$ almost everywhere. In
particular, $t(\lunl g\runl ,W)=0$.
\end{lemma}

The condition that $t_k(g,W)=0$ for some $k$-labeled quantum graph
$g$ (where $k>0$) seems to carry much more information than a
condition that $t(f,W)=0$ for an unlabeled quantum graph $f$.
However, there is a way to translate labeled constraints to unlabeled
constraints. First we state a simple version.

\begin{lemma}\label{LEM:UNLABEL1}
Let $f$ be a $k$-labeled quantum graph and $d$, a positive even
integer. Let $f^d$ be the $d$-th power of $f$ in the algebra $\GG_k$.
Then for any $W\in\WW$, $t(\lunl f^d\runl ,W)=0$ if and only if
$t_k(f,W)=0$ almost everywhere.
\end{lemma}

\noindent(By Lemma \ref{LEM:UNLABEL0}, the same conclusion holds if
only some of the nodes are unlabeled.)

\begin{proof}
We have
\[
t(\lunl f^d\runl ,W)=\int\limits_{[0,1]^k} t_k(f^d,W)(x)\,dx =
\int\limits_{[0,1]^k}\bigl( t_k(f,W)(x)\bigr)^d\,dx,
\]
so this is $0$ if and only if $t_k(f,W)(x)=0$ for almost all
$x\in[0,1]^k$.
\end{proof}

Let $F$ be a $k$-labeled multigraph, and let $\PP=\{S_1,\dots,S_m\}$
be a partition of $[k]$. We say that $\PP$ is {\it legitimate} for
$F$, if each set $S_i$ is independent in $F$. If this is the case,
then we define the $m$-labeled multigraph $F/\PP$ by identifying the
nodes in each $S_i$, and labeling the obtained node with $i$. For a
$k$-labeled quantum graph $f$, we say that the partition $\PP$ of
$[k]$ is {\it legitimate for $f$} if it is legitimate for every
constituent. Then we can define $f/\PP$ by linear extension.

\begin{lemma}\label{LEM:CONTRACT}
Let $f$ be a $k$-labeled quantum graph and $\PP$, a legitimate
partition for $f$ with $r$ classes. Let $W\in\WW$, and suppose that
$t_k(f,W)=0$ almost everywhere. Then $t_r(f/\PP,W)=0$ almost
everywhere.
\end{lemma}

\begin{proof}
If $k=2$ and $\PP$ identifies the two labels, then the Lemma follows
from Lemma \ref{derand2}, since $t_1(f/\PP,W)(x)=t_2(f,W)(x,x)$ is
$0$ whenever $x$ is a typical point, i.e., almost everywhere. (We
could also invoke Theorem 1.6 in \cite{LSz3} here.)

Now for an arbitrary $k\ge 2$, it suffices to prove the case when
$\PP$ identifies a single pair of labeled nodes, say $1$ and $2$ (so
$r=k-1$). If $t_k(f,W)=0$ almost everywhere, then by
\eqref{EQ:MULTIPL} $t_k(f^2,W)=0$ almost everywhere. Let $h$ be
obtained from $f^2$ by unlabeling $3,\dots,k$, and let $h'$ denote
the $1$-labeled quantum graph obtained from $h$ by identifying labels
$1$ and $2$. Then by Lemma \ref{LEM:UNLABEL0} it follows that
$t_2(h,W)=0$ almost everywhere. By the above, we have $t_1(h',W)=0$
almost everywhere. Again by Lemma \ref{LEM:UNLABEL0}, we have
$t(\lunl h'\runl ,W)=0$. But $\lunl h'\runl =\lunl (f/\PP)^2\runl $,
and hence by Lemma \ref{LEM:UNLABEL1} it follows that
$t_{k-1}(f/\PP,W)=0$ almost everywhere.
\end{proof}

One drawback of Lemma \ref{LEM:UNLABEL1} is that $f^r$ may have
multiple edges, even if $f$ does not. The construction in the next
Lemma gets around this.

\begin{lemma}\label{LEM:UNLABEL2}
For every simple $k$-labeled quantum graph $f$ there is a simple
unlabeled quantum graph $g$ such that for any $W\in\WW$, $t(g,W)=0$
if and only if $t_k(f,W)=0$ almost everywhere.
\end{lemma}

\begin{proof}
For every $k$-labeled quantum graph $g$, consider the product (in the
algebra $\GG_k$) of all constituents, and define $\Lab(g)$ as the
subgraph of this induced by the labeled nodes.

We prove the lemma by induction on the chromatic number
$\chi(\Lab(f))$. If $\chi(\Lab(f))=1$, then the labeled nodes are
independent in every constituent, and hence we can take $g=\lunl
f^2\runl $ and use Lemma \ref{LEM:UNLABEL1}.

Suppose that  $\chi(\Lab(f))=r>1$, and let $[k]= S_1\cup\dots\cup
S_r$ be an $r$-coloring of $\Lab(f)$. We can assume that
$S_r=\{k-q+1,\dots,k\}$. We glue together two copies of $f$ along
$S_r$. Formally, let $f_1$ be obtained from $f$ by increasing the
labels in $S_r$ by $k-q$ (the labels not in $S_r$ are not changed),
and by adding isolated nodes labeled $k-q+1,\dots,2k-2q$ to every
constituent of $f$. Let $f_2$ be obtained from $f$ by increasing all
labels by $k-q$, and by adding isolated nodes labeled $1,\dots,k-q$
to every constituent of $f$. So $f_1$ and $f_2$ are $(2k-q)$-labeled
quantum graphs. Form the product $f_1f_2$ and remove the labels
$2k-2q+1,\dots,2k-q$, to get a $(2k-2q)$-labeled quantum graph $h$.

\begin{claim}\label{CLAIM:FG}
For every $W\in\WW$, $t_k(f,W)=0$ almost everywhere if and only if
$t_{2k-2q}(h,W)=0$ almost everywhere.
\end{claim}

The ``only if'' part is obvious, since $t_k(f,W)=0$ almost everywhere
implies $t_{2k-q}(f_1,W)=t_{2k-q}(f_2,W)=0$ almost everywhere, which
implies $t_{2k-q}(f_1f_2,W)=0$ almost everywhere, which in turn
implies $t_{2k-2q}(h,W)=0$ almost everywhere.

To prove the ``if'' part, note that two labeled nodes whose labels
correspond to the same label in $f$ are never adjacent, so we can
identify these labels in $h$ to get $f^2$ (with the labels in $S_r$
removed). So $t_{2k-2q}(h,W)=0$ almost everywhere implies by Lemmas
\ref{LEM:CONTRACT} and \ref{LEM:UNLABEL0} that $t(\lunl f^2\runl
,W)=0$, and hence by Lemma \ref{LEM:UNLABEL1}, we get that
$t_k(f,W)=0$ almost everywhere. This proves the Claim.

Thus it suffices to replace the constraint $t_{2k-2q}(h,W)=0$ by an
unlabeled constraint. This can be done by induction, since
$\chi(\Lab(h))\le r-1$.
\end{proof}

\begin{corollary}\label{COR:LABFORCE}
Suppose that for $W\in\WW$ there is a family $\{f_1,\dots,f_m\}$,
where $f_i$ is a simple $k_i$-labeled quantum graph such that
$t_{k_i}(f_i,W)=0$ almost everywhere, and the constraints
$t_{k_i}(f_i,U)=0$, $U\in\WW$ imply that $U$ is weakly isomorphic to
$W$. Then $W$ is finitely forcible in $\WW$. Similar assertion holds
for forcing in $\WW_0$.
\end{corollary}

\subsection{The adjoint of an operator}

Let $\FF$ denote the set of simple graphs (up to isomorphism), and
let $\QQ$ be the linear space of simple quantum graphs.

\begin{definition}
Let $\mathbf{F}:~\WW \rightarrow\WW $ be an operator (not necessarily
linear) preserving weak isomorphism, and let
$\mathbf{F}^*:~\QQ\rightarrow \QQ$ be a linear map. We say that the
map $\mathbf{F}^*$ is an {\it adjoint} of $\mathbf{F}$ if
\[
t(g,\mathbf{F}(W))=t(\mathbf{F}^*(g),W)
\]
for every $g\in\mathcal{Q}$ and $W\in\WW $. (Note that it is enough
to define $\mathbf{F}^*$ on simple graphs and extend it linearly to
quantum graphs.) We denote the set of functionals which have an
adjoint by $\DD$.
\end{definition}

It is clear from this definition that the elements of $\DD$ form a
semigroup with respect to composition.

\begin{example}\label{DUAL-MULTI}
Fix a real number $\alpha$ and let $\mathbf{F}$ denote the functional
defined by $\mathbf{F}(W)=\alpha W$. It is easy to see that
$\mathbf{F}$ has an adjoint defined for simple graphs $G$ by
$\mathbf{F}^*(G)=\alpha^{|E(G)|}G$.
\end{example}

\begin{example}\label{DUAL-ADD}
Let $\beta$ be a real number and $\mathbf{F}(W)=W+\beta$. Then
$\mathbf{F}$ has an adjoint defined by
\[
\mathbf{F}^*(G)=\sum_{Z\subseteq E(G)} \beta^{|E(G)\setminus Z|}
(V(G),Z).
\]
\end{example}

\begin{example}
Let $U\in\WW $ be a fixed function and define $\mathbf{F}(W)$ as the
tensor product $U\otimes W$. Then $\mathbf{F}^*(G)=t(G,U)G$ defines
an adjoint of $\mathbf{F}$.
\end{example}

\begin{example}
Let $\mathbf{F}(W)=W^{\otimes k}$ be the $k$-th tensor power of $W$
(for a fixed $k\ge 1$). Then an adjoint $\mathbf{F}^*$ can be defined
by letting $\mathbf{F}^*(G)$ be the disjoint union of $k$ copies of
$G$.
\end{example}

\begin{example}
Let $p(z)=\sum_{k=1}^n a_k z^k$ be a real valued polynomial. We
define $\mathbf{F}(W)$ as $p(W)$ where $W$ is substituted into $p$ as
an integral kernel operator. For any graph $G=(V,E)$, we define
\[
\mathbf{F}^*(G)=\sum_{\mathbf{k}\in[n]^E} a_\mathbf{k}
G^{(\mathbf{k})},
\]
where for $\mathbf{k}\in [n]^E$ we define $a_\mathbf{k}=\prod_{e\in
E} a_{k_e}$, and $G^{(\mathbf{k})}$ is the graph obtained from $G$ by
subdividing each edge $e$ by $k_e-1$ nodes.

We show that $\mathbf{F}^*$ is an adjoint of $\mathbf{F}$. We use
that for $\mathbf{k}\in[n]^E$
\[
t(G^{(\mathbf{k})},W) = \int\limits_{[0,1]^{V(G)}} \prod_{ij\in E(G)}
W^{\circ k_{ij}}(x_i,x_j)\,dx.
\]
Hence
\begin{align*}
t(G,p(W))&=\int\limits_{[0,1]^{V(G)}} \prod_{ij\in E(G)}
\Bigl(\sum_{k=1}^n a_k W^{\circ k}(x_i,x_j)\Bigr)\,dx\\
&= \int\limits_{[0,1]^{V(G)}} \sum_{\mathbf{k}\in[n]^E} a_\mathbf{k}
\prod_{ij\in E(G)} W^{\circ k_{ij}}(x_i,x_j)\,dx=
\sum_{\mathbf{k}\in[n]^E} t(G^{(\mathbf{k})},W) =
t(\mathbf{F}^*(G),W).
\end{align*}
\end{example}

\begin{example}
Let $H$ be a simple 2-labeled graph which has an automorphism
interchanging the labeled nodes. Then
$\mathbf{F}_H(W)(x,y)=t_2(H,W)(x,y)$ is a symmetric 2-variable
function in $x,y$. Let $\mathbf{F}_H^*(G)$ be the graph obtained from
$G$ by replacing each edge by a copy of $H$ where the labeled nodes
of $H$ are identified with the endpoints of the edge. Then
$\mathbf{F}_H^*$ is an adjoint of $\mathbf{F}_H$.

As a special case, if $H=K_3^{\bullet\bullet}$ denotes the triangle
with two labeled nodes, then $\mathbf{F}_H(W)=(W\circ W)W$, and
$\mathbf{F}_H^*(G)$ can be constructed by doubling each edge of $G$
and subdividing one copy of each edge.
\end{example}

\begin{lemma}\label{fininv}
Let $W\in\WW$ be finitely forcible in $\WW$, let $\mathbf{F}\in\DD$,
and assume that $\mathbf{F}^{-1}([W])$ is finite (up to weak
isomorphism). Then every element in $\mathbf{F}^{-1}([W])$ is
finitely forcible in $\WW$.
\end{lemma}

\begin{proof}
If $W$ can be forced by the constraints $t(F_i,W)=a_i$
$(i=1,\dots,k)$, then the set $\mathbf{F}^{-1}([W])$ can be forced by
the constraints $t(\mathbf{F}^*(F_i),U)=a_i$. Let
$\mathbf{F}^{-1}([W])$ have $m$ elements up to weak isomorphism, then
the equivalence class of each element can be distinguished from the
others by at most $m-1$ further graph density constraints.
\end{proof}

Applying this lemma with examples \ref{DUAL-MULTI} and
\ref{DUAL-ADD}, we get

\begin{corollary}
If $W\in\WW $ is finitely forcible (in $\WW$), then so is $\alpha
W+\beta$ for $\alpha,\beta\in\mathbb{R}$.
\end{corollary}

\begin{corollary}\label{lintr}
For every finitely forcible graphon $W\in\WW$ there are numbers
$\alpha\neq 0$ and $\beta$ such that $\alpha W+\beta$ is in $\WW_0$
and is finitely forcible in $\WW_0$.
\end{corollary}

For our next corollary, we need a simple lemma. (Here we substitute a
graphon into a polynomial as a kernel operator.)

\begin{lemma}\label{BIJECT}
Let $p$ be a polynomial which is a bijection on $\mathbb{R}$ with
$p(0)=0$. Then

\smallskip

{\rm(a)} If $p(W_1)=p(W_2)$ almost everywhere, then $W_1=W_2$ almost
everywhere.

\smallskip

{\rm(b)} If $p(W_1)$ and $p(W_2)$ are weakly isomorphic, then so are
$W_1$ and $W_2$.
\end{lemma}

\begin{proof}
(a) Let $U\in\WW $, and consider any function $W\in\WW$ with $p(W)=U$
almost everywhere. Let
\[
W(x,y)\sim \sum_{i=1}^\infty \mu_i f_i(x)f_i(y).
\]
be the spectral decomposition of $W$, where $\{f_i\}_{i=1}^\infty$ is
an orthonormal system of functions in $L_2[0,1]$ and $\sum_i \mu_i^2
<\infty$. Then
\[
U(x,y)\sim \sum_{i=1}^\infty p(\mu_i) f_i(x)f_i(y),
\]
is a spectral decomposition of $U$, since $p$ is Lipshitz on any bounded
interval and hence $\sum_i p(\mu_i)^2
<\infty$. Since the spectral decomposition of $U$ is unique (up to an
orthogonal basis transformation in the eigensubspaces) and $p$ is
injective, and we see that the $\mu_i$ and $f_i$ are determined by
$U$ (again, up to an orthogonal basis transformation in the
eigensubspaces), and so $W$ is determined by $U$.

(b) Assume that $p(W_1)$ and $p(W_2)$ are weakly isomorphic. By
Theorem \ref{THM:WEAKISO}, this implies that there are measure
preserving maps $\varphi,\psi:~[0,1]\to[0,1]$ such that
$p(W_1)^\varphi = p(W_2)^\psi$ almost everywhere. It is easy to check
that $p(W_1)^\varphi = p(W_1^\varphi)$, so we get that
$p(W_1^\varphi) = p(W_2^\psi)$ almost everywhere. By (a), this means
that $W_1^\varphi = W_2^\psi$ almost everywhere, and so $W_1$ and
$W_2$ are weakly isomorphic.
\end{proof}

\begin{corollary}\label{invforc}
Let $p$ be a polynomial which is a bijection on $\mathbb{R}$ with
$p(0)=0$. If $p(W)$ is finitely forcible for some $W\in\WW $, then so
is $W$.
\end{corollary}

\begin{proof}
Let $\mathbf{F}(W)=p(W)$. By Lemma \ref{BIJECT}(b),
$\mathbf{F}^{-1}([p(W)])$ is finite up to weak isomorphism (in fact, has
at most one element). Hence by Lemma \ref{fininv}, $W$ is finitely
forcible.
\end{proof}

\section{Finitely forcible properties}\label{PROPS}

As mentioned in the introduction, instead of forcing specific
graphons by a finite number of subgraph densities, we can more
generally ask which properties of graphons can be forced this way.
Clearly, every such property is invariant under weak isomorphism, and
also closed under convergence. (More generally, it is closed in the
cut-norm \cite{LSz1,BCLSV1}, but we don't need this in this paper.)

Some important properties are finitely forcible, but some others are
not. It is sometimes the case, however, that in the presence of some
other condition, such properties become finitely forcible. The
property that $W$ is $0/1$ valued is an example (to be
discussed below).

\subsection{Regularity}

We call a graphon {\it $d$-regular}, or {\it regular of degree $d$}
($0\le d\le 1$), if
\[
\int\limits_0^1 W(x,y)\,dy=d
\]
for almost all $0\le x\le 1$. This condition can also be written as
$t_1(K_2^\bullet,W)(x)=d$, where $K_2^\bullet$ denotes the single
edge with one endnode labeled. These graphons can be forced by two
subgraph density constraints:
\[
t(K_2,W)=d,\qquad t(P_3,W)=d^2,
\]
since equality holds in the Cauchy--Schwarz estimate
\[
t(P_3,W) = \int\limits_0^1 t(K_2^\bullet,W)^2 \,dx \ge \Bigl(
 \int\limits_0^1 t(K_2^\bullet,W)\,dx\Bigr)^2 = t(K_2,W)^2.
\]
Regular graphons (without specifying the degree $d$) can be forced by
the constraint $t(P_3,W)=t(K_2,W)^2$.

\subsection{Zero-one valued functions}

Trivially, $W\in\WW$ is $0/1$ valued almost everywhere if and
only if $t_2(\widehat{C}_2,W)=0$, where $\widehat{C}_2$ is the
2-labeled signed multigraph on 2 nodes with 2 parallel edges, one
signed ``$+$'' and one signed ``$-$''. By Lemma \ref{LEM:UNLABEL1},
this is equivalent to the single numerical equation
$t(\widehat{B}_4,W)=0$, where $\widehat{B}_4$ is the unlabeled signed
multigraph on 2 nodes with 4 parallel edges, 2 signed ``$+$'' and 2
signed ``$-$''.

So we can ``force'' the property of being $0/1$ valued using
multigraphs, but we cannot express it in terms of simple graphs. This
follows from the observation that if $G(n,1/2)$ is the
Erd\H{o}s--R\'enyi random graph with $n$ nodes and edge density
$1/2$, and $W_n=W_{G(n,1/2)}$, then with probability $1$, $W_n$ tends
to the identically $1/2$ function $U_{1/2}$ in the $\|.\|_\square$
norm, which implies by \eqref{EQ:COUNT} that $t(F,W_n)\to
t(F,U_{1/2})$ for every simple graph $F$. So every constraint of the
form $t(F,W_n)=0$, where $F$ is a simple graph, is inherited by
$U_{1/2}$, which is not $0/1$ valued.

It makes sense to formulate sufficient conditions for being
$0/1$ valued. Here is a useful one.

\begin{lemma}\label{M2FREE}
Let $\widehat{F}$ be a signed bipartite graph on $n$ nodes, all
labeled. Let $\lunl \widehat{F}\runl $ be obtained by unlabeling all
the nodes. Suppose that for some $W\in\WW$ we have
\begin{equation}\label{MON-1}
t_n(\widehat{F},W)=0
\end{equation}
almost everywhere. Then $W(x,y)\in\{0,1\}$ almost everywhere. If
$W\in\WW_0$, then it suffices to assume that
\begin{equation}\label{MON-2}
t(\lunl \widehat{F}\runl ,W)=0.
\end{equation}
\end{lemma}

\begin{proof}
By Lemma \ref{LEM:CONTRACT}, \eqref{MON-1} implies that for the
2-labeled signed multigraph $J$ obtained by identifying each color
class of $F$, we have $t_2(J,W)=0$. This clearly implies that $W$ is
$0/1$ valued.

If $W\in\WW_0$, then \eqref{MON-2} implies that in the integral
\[
t(\lunl \widehat{F}\runl
,W)=\int\limits_{[0,1]^{V(F)}}t_n(\widehat{F},W)(x)\,dx
\]
the integrand is $0$ almost everywhere. This means that \eqref{MON-1}
holds.
\end{proof}

\subsection{Monotonicity}

Let $\MM_0$ denote the set of measurable functions
$[0,1]^2\to\{0,1\}$ that are monotone decreasing in both variables,
and let $\MM$ be the set of graphons which are weakly isomorphic to
some function in $\MM_0$. (These graphons have been studied by
Diaconis, Holmes and Janson \cite{DHJ} as limits of threshold
graphs.) In this section we show that the set $\mathcal{M}$ is
finitely forcible in $\WW$.

Let $\widehat{C}_4$ denote a signed 4-labeled 4-cycle, with two
opposite edges signed ``$+$'', the other two signed ``$-$''. Let
$\lunl \widehat{C}_4\runl $ be obtained from $\widehat{C}_4$ by
unlabeling all its nodes.

\begin{lemma}\label{M2FREE2}
Let $W\in\WW$, then $W\in\MM$ if and only if
\begin{equation}\label{EQ:T4C4}
t_4(\widehat{C_4},W)=0
\end{equation}
almost everywhere. If $W\in\WW_0$, then it is enough to assume that
\begin{equation}\label{EQ:TC4}
t(\lunl \widehat{C_4}\runl ,W)=0
\end{equation}
almost everywhere.
\end{lemma}

\begin{proof}
It is easy to see that \eqref{EQ:T4C4} and \eqref{EQ:TC4} hold for
every $W\in\MM$.

Next we prove if that $W\in\WW_0$ satisfies \eqref{EQ:TC4} then
$W\in\MM$. By deleting zero-sets, we may assume that all points are
typical. By Lemma \ref{M2FREE}, $W$ is $0/1$ valued almost
everywhere, and we may assume that it is $0/1$ valued. Let $N(x)$
denote the support of the function $W(x,.)$ ($x\in[0,1]$).

By the Monotone Reordering Theorem, there is a monotone decreasing
function $f:~[0,1]\to[0,1]$ and a measure preserving map
$\varphi:~[0,1]\to[0,1]$ such that $\lambda(N(x))=f(\varphi(x))$
almost everywhere. We can change $W$ on a set of measure $0$ so that
$\lambda(N(x))=f(\varphi(x))$ for all $x$. Consider the function
$U(x,y)=\one(y\le f(x),x\le f(y))$. This is clearly symmetric and
monotone decreasing in both variables. We claim that $W=U^\varphi$
almost everywhere.

\begin{claim}\label{CLAIM:ORDERED}
For all $x,y\in[0,1]$, either $\lambda(N(x)\setminus  N(y))=0$ or
$\lambda(N(y)\setminus  N(x))=0$.
\end{claim}

Indeed, we have
\[
\int W(x,u)(1-W(u,y))W(y,v)(1-W(v,x))\,dx\,dy\,dz\,du = t(\lunl
\widehat{C_4}\runl ,W)=0,
\]
and hence by the fact that $W$ is $0$-$1$ valued, we have
$W(x,u)(1-W(u,y))W(y,v)(1-W(v,x))=0$ almost everywhere. Hence for
almost every pair $(x,y)$, we have $\int
W(x,u)(1-W(u,y))W(y,v)(1-W(v,x))\,du\,dv=0$. By Lemma \ref{derand2},
this implies that this holds for all $x$ and $y$. But this integrand
is positive if $u\in N(x)\setminus  N(y)$ and $v\in N(y)\setminus
N(x)$, and so one of these sets must have measure $0$.

Let
\[
A= \{(x,y,u)\in[0,1]^3:~\varphi(x)\le \varphi(y), ~u\in N(y)\setminus
N(x)\}
\]
Since $f$ is monotone decreasing, $\varphi(x)\le\varphi(y)$ implies
that $\lambda(N(x))\ge \lambda(N(y))$, and so Claim
\ref{CLAIM:ORDERED} implies that $\lambda(N(y)\setminus N(x))=0$.
Hence $\lambda\{u:~(x,y,u)\in A\}=0$ for almost all pairs $(x,y)$,
which implies that $\lambda_3(A)=0$.

Let $K(u)=\varphi^{-1}[0,f(\varphi(u))]$. We have $\lambda(N(u))
=f(\varphi(u))=\lambda(K(u))$ for every $u$. We claim that
$\lambda(N(u)\triangle K(u))=0$ for almost all $u$. Indeed, if
$\lambda(N(u)\triangle K(u))>0$, then $\lambda(N(u)\setminus
K(u))=\lambda(K(u)\setminus N(u))>0$, and then for every $x\in
K(u)\setminus N(u)$ and $y\in N(u)\setminus K(u)$ we have $\phi(x)\le
f(\phi(u))<\phi(y)$ and $u\in N(y)\setminus N(x)$, and so $(x,y,u)\in
A$. Since $\lambda_3(A)=0$, this can hold for a nullset of points $u$
only.

So we know that for almost all $y$, $\lambda(N(y)\triangle K(y))=0$.
Hence for almost all pairs $(x,y)$, $W(x,y)=1$ if and only if $x\in
K(y)$, i.e., $\varphi(x)\le f(\varphi(y))$. By the symmetry of $W$,
this is also equivalent (for almost all pairs $x,y$) to
$\varphi(y)\le f(\varphi(x))$. So (for almost all pairs $x,y$)
$W(x,y)=1$ if and only if $U(\varphi(x),\varphi(y))=1$. This shows
that $W$ is weakly isomorphic to $U\in\MM_0$.

Finally, assume that $W\in\WW$ satisfies \eqref{EQ:T4C4}. By Lemma
\ref{M2FREE}, $W$ is $0/1$ valued almost everywhere, so
$W\in\WW_0$. Since it trivially satisfies \eqref{EQ:TC4}, it follows
that $W\in\MM$.
\end{proof}

\section{Finitely forcible graphons I: polynomials}

\subsection{Positive supports of polynomials}\label{POLYN}

\begin{theorem}\label{MON-POLY}
Let $p$ be a real symmetric polynomial in two variables, which is
monotone decreasing on $[0,1]^2$. Then the function $W(x,y)
=\one_{p(x,y)\ge 0}$ is finitely forcible in $\WW$.
\end{theorem}

\begin{proof}
We in fact prove that the equations
\begin{equation}\label{C4W}
t_4(\widehat{C_4},U) =0
\end{equation}
and
\begin{equation}\label{KABW}
t(K_{a,b},U) = t(K_{a,b},W)\qquad (1\le a,b\le 2\deg(p)+2)
\end{equation}
form a forcing family for $W$ in $\WW$. The condition on the
monotonicity of $p$ implies that $U=W$ satisfies \eqref{C4W}. It is
trivial that the other equations are satisfied by $U=W$.

Let $U\in\WW$ be any graphon satisfying \eqref{C4W}-\eqref{KABW}. By
Lemma \ref{M2FREE2} we may assume that $U$ is $0/1$ valued and
monotone decreasing. Let $S_U=\{(x,y):~U(x,y)=1\}$.

We have
\[
t(K_{a,b},U)= \int\limits_{[0,1]^a}\int\limits_{[0,1]^b}
\prod_{i=1}^a \prod_{j=1}^b U(x_i,y_j)\,dy\,dx.
\]
Split this integral according to which $x_i$ and which $y_j$ is the
largest. Restricting the integral to, say, the domain where $x_1$ and
$y_1$ are the largest, we have that whenever $U(x_1,y_1)=1$ then also
$U(x_i,y_j)=1$ for all $i$ and $j$, and hence
\begin{align*}
\int\limits_{x_1\in [0,1]} & \int\limits_{~x_2,\dots,x_a\le x_1}
\int\limits_{~y_1\in [0,1]} \int\limits_{~y_2,\dots,y_b \le y_1}
\prod_{i=1}^a \prod_{j=1}^b
U(x_i,y_j)\,dy\,dx\\
&= \int\limits_{x_1\in [0,1]} \int\limits_{~y_1\in [0,1]} U(x_1,y_1)
x_1^{a-1} y_1^{b-1} \,dy_1\,dx_1= \int\limits_{(x,y)\in S_U} x^{a-1}
y^{b-1} \,dy\,dx\,.
\end{align*}
Hence
\[
t(K_{a,b},U)=ab \int\limits_{(x,y)\in S_U} x^{a-1} y^{b-1}
\,dy\,dx\,.
\]
By Stokes' Theorem, we can rewrite this as
\[
t(K_{a,b},U)= b \int\limits_{\partial S_U} x^a y^{b-1} n_1(x,y)
\,ds\,,
\]
where $ds$ is the arc length of $\partial S_U$ and $n=(n_1,n_2)$ is
the outward normal of $\partial S_U$. (Since $\partial S_U$ is the
graph of a monotone function, this normal exists almost everywhere.)
Interchanging the roles of $x$ and $y$, and adding, we get
\begin{equation}\label{MONOM}
\int\limits_{\partial S_U} x^a y^b (n_1(x,y)+n_2(x,y)) \,ds
=\frac{1}{a+1} t(K_{a+1,b},U)+\frac{1}{b+1} t(K_{a,b+1},U).
\end{equation}

Now consider the following integral:
\[
I(U)=\int\limits_{\partial S_U} xyp(x,y)^2 (n_1(x,y)+n_2(x,y)) \,ds.
\]
By \eqref{MONOM}, this can be expressed as a linear combination of
the values $t(K_{a,b},U)$, where $a,b\le 2\deg(p)+1$ and the
coefficients depend only on $a,b$ and $p$. Hence it follows that
$I(U)=I(W)=0$.

On the other hand, the integrand in $I(U)$ is clearly $0$ on the
axes, and it is nonnegative on the rest of the boundary. Hence it
must be identically $0$, which means that $\partial(S_U)$ must be
contained in the union of the axes and the curve $p=0$. But this
clearly implies that $U=W$ except perhaps on the boundary.
\end{proof}

The following special case is perhaps the simplest. Define the {\it
half-graphon} by $W_h(x,y)=\one_{x+y\le 1}$.

\begin{corollary}\label{haromszog}
The half-graphon is finitely forcible in $\WW$.
\end{corollary}

In fact, by a variation of the argument above, one can prove that the
following equations force the half-graphon:
\begin{equation}\label{HALF1}
t(\widehat{C}_4,W)=0,
\end{equation}
and
\begin{equation}\label{HALF2}
t(P_3,W)-t(K_2,W)+1/6=0.
\end{equation}
Clearly, the left hand side of \eqref{HALF1} is always nonnegative.
It is easy to show that in \eqref{HALF2}, the left had side is
nonnegative, provided equality holds in \eqref{HALF1}.

Half-graphons are natural limits of {\it half-graphs}, defined by
$V(G)=[n]$ and $E(G)=\{ij:~i+j\le n\}$. This
implies the following graph-theoretic extremal result.

\begin{corollary}\label{HALFGR}
Among all simple graphs with no induced matching with $2$ edges, the
difference $t(P_3,.)-t(K_2,.)$ is asymptotically minimal for
half-graphs.
\end{corollary}

\subsection{Continuous range}

Applying the results of Section \ref{OPER} to the half-graphon
$U(x,y)=\one_{x+y>1}$, we get an interesting finitely forcible
graphon.

\begin{prop}\label{CONTIN}
There exists a graphon $W\in\WW $ satisfying $W+\frac12(W\circ W\circ
W)=U$, it is finitely forcible in $\WW$, and its range consists of
two nontrivial intervals.
\end{prop}

\begin{proof}
Let $\lambda$ be a non-zero eigenvalue of $U$ and let $f$ be the
corresponding eigenfunction with unit norm. Then
\[
\lambda f(x)=\int\limits_{1-x}^1 f(y)~dy.
\]
From this integral equation one can calculate that the eigenvalues
and eigenfunctions are
\[
\lambda_k=\frac{2}{(4k+1)\pi}, \qquad
f_k(x)=\sqrt{2}\sin\Bigl(\frac{4k+1}{2}\pi x\Bigr),\qquad k\in\Z.
\]
In particular, the eigenfunctions are analytic and uniformly bounded.

Let $g$ be the inverse function of $x\mapsto x+\frac12x^3$ on the
real line. We have
\[
x-g(x) = \frac{xg(x)^2}{2+ g(x)^2},
\]
which implies that if $x\neq 0$ then $|g(x)|<|x|$ and $|x-g(x)|<\frac12 |x|
g(x)^2<\frac12 |x|^3$. Since the series $\sum_k |\lambda_k|^3$ is
convergent, and the $f_k$ are bounded, this implies that
\[
P(x,y)=\sum_{k\in\Z}(\lambda_k-g(\lambda_k))f_k(x)f_k(y)
\]
is uniformly absolute convergent. Thus $P(x,y)$ is well defined and
analytic in $x$ and $y$. It is clear that $P$ is not constant.

Note that we have
\begin{equation}\label{EQ:PA}
|P(x,y)|\le \sum_{k\in\Z} 2 |\lambda_k-g(\lambda_k)|< \sum_{k\in\Z}
|\lambda_k|^3 = \Bigl(\frac2{\pi}\Bigr)^3\frac78\zeta(3)< \frac12.
\end{equation}

Now we can define $W=U-P$. Clearly this function is bounded and
symmetric. Furthermore, the spectral decomposition of $W$ is
\[
W(x,y)=\sum_{k\in\Z} \lambda_kf_k(x)f_k(y)-\sum_{k\in\Z}
(\lambda_k-g(\lambda_k))f_k(x)f_k(y) =\sum_{k\in\Z}
g(\lambda_k)f_k(x)f_k(y),
\]
from which it follows that $W+\frac12(W\circ W\circ W)=U$. Corollary
\ref{invforc} implies that $W$ is finitely forcible.

Clearly $W$ is continuous over $U^{-1}(0)$ and over $U^{-1}(1)$.
Inequality \eqref{EQ:PA} implies that the range of $W$ over
$U^{-1}(0)$ is an interval contained in $(-1/2,1/2)$, while its range
over $U^{-1}(1)$ is an interval contained in $(1/2,3/2)$. So we get
two disjoint intervals.
\end{proof}

We can invoke Corollary \ref{lintr} to transform $W$ into an element
from $\WW_0$. This implies the following.

\begin{corollary}
There is finitely forcible function in $\WW_0$ whose range is of
continuum cardinality.
\end{corollary}

\section{Finitely forcible graphons II: complement reducible graphons}

\subsection{The finite case}

A simple graph is called {\it complement reducible}, or for short, a
{\it CR-graph}, if it can be constructed starting from a single node,
by repeated application of disjoint union and complementation. (These
graphs are often called {\it cographs}). One of the many known
characterizations of these graphs is the following (\cite{CLB}; see
\cite{BLS} for more on these graphs).

\begin{prop}\label{PROP:CR}
A simple graph is complement reducible if and only if it does not
contain a path $P_4$ on $4$ nodes as an induced subgraph.
\end{prop}

Let $\widehat P_4$ denote the graph $K_4$ in which the edges of a
path of length $3$ are signed ``$+$'', and the remaining edges are
signed ``$-$''. Then the condition in the proposition can be
rephrased as
\[
t(\widehat{P}_4,G)=0.
\]

Every CR-graph $G$ can be described by a rooted tree $T$ in a natural
way \cite{CLB}: each node of $T$ represents a CR-graph; the leaves
represent single nodes, the root represents $G$, and the children
of each internal node represent the connected components of the graph
represented by the node, complemented. Another way of describing this
connection is that $G$ is defined on the leaves of $T$, and two nodes
of $G$ are connected by an edge if and only if their last common
ancestor is at an odd distance from the root.

\subsection{CR graphons and trees}

We define a {\it CR-graphon} as any graphon $W$ with
$t(\widehat{P}_4,W)=0$. Also we extend the notion of CR-graphs to
infinite graphs by the requirement that no four nodes induce a path
of length four.

One can construct CR-graphons from trees, but (unlike in the finite
case) not all CR-graphons arise this way (we shall see later that all
regular CR-graphons do). Let $T$ be any (possibly infinite) rooted
tree, in which every non-leaf node has at least two and at most
countably many children, except possibly the root, which may have
only one child. Let $\Omega=\Omega_T$ be the set of maximal paths
starting at the root $r$ (they are either infinite or end at a leaf).
For each node $v$, let $C_v$ be the set of its children, and let
$\Omega_v$ denote the set of paths in $\Omega$ passing through $v$.
The sets $\Omega_v$ generate a $\sigma$-algebra $\AA_T$.

We define a (simple) graph on node set $\Omega$ by connecting two
nodes if the last common node of the corresponding paths is at odd
depth (where the depth of the root $r$ is $0$). We also define the
adjacency function $U_T:~\Omega\times \Omega\to\{0,1\}$ by letting
$U_T(x,y)=1$ if and only if $x$ and $y$ are adjacent. It is clear
that $U_T$ is measurable with respect to $\AA\times\AA$.

If we choose any probability measure $\pi$ on $(\Omega_T,\AA_T)$,
this completes the construction of a graphon
$(\Omega_T,\AA_T,\pi,U_T)$, which is clearly a CR-graphon. Note that
such a measure can be specified through the values
\[
f(v)=\pi(\Omega_v).
\]
It is clear that these values satisfy
\begin{equation}\label{EQ:F-ADD}
f(r)=1, \qquad f(u)\ge 0, \qquad\text{and}\qquad f(u)=\sum_{v\in C_u}
f(v).
\end{equation}
Conversely, every function satisfying \eqref{EQ:F-ADD} defines a
probability measure on $(\Omega_T,\AA_T)$.

\subsection{Regular CR-graphons}

Of special interest for us will be regular CR-graphons. Our first
goal is to prove:

\begin{theorem}\label{THM:REG-LOCFIN}
Every regular CR-graphon $W$ can be represented (up to weak
isomorphism) as $(\Omega_T,\AA_T,\pi,U_T)$ where $T$ is a locally
finite tree and $\pi$ is a probability measure on $\AA_T$.
\end{theorem}

The proof of this theorem will need several lemmas. Recall from
\cite{LSz1} that for every graphon $W:~[0,1]^2\mapsto [0,1]$ there is
a random graph model $\mathbb{G}(W,n)$ on node set $\{1,2,\dots,n\}$,
created as follows: We pick independent uniform random points
$x_1,x_2,\dots,x_n\in[0,1]$ and connect two distinct nodes $i$ and
$j$ with probability $W(x_i,x_j)$. Let $d_i(\mathbb{G}(W,n))$ denote
the degree of $i$ in the resulting graph.

\begin{lemma}\label{degcont}
If a graphon $W$ is $d$-regular, then with probability at least
$1-2ne^{-(n-1)\eps^2/2}$ we have
\[
\Bigl| \frac{d_i(\mathbb{G}(W,n))}{n-1}-d \Bigr|<\eps
\]
simultaneously for all $1\leq i\leq n$.
\end{lemma}

\begin{proof}
It is easy to see that for every $1\leq i\leq n$ the value
$d_i(\mathbb{G}(W,n))$ is the sum of $n-1$ independent random
variables all taking $1$ with probability $d$ and $0$ with
probability $1-d$. This implies by the Chernoff--Hoeffding Inequality
that
\[
P\Bigl(\Bigl|\frac{d_i(\mathbb{G}(W,n))}{n-1}-d\Bigr|\geq\eps\Bigr)\leq
2e^{-(n-1)\eps^2/2}.
\]
This means that the probability that there exist at least one number
$1\leq i\leq n$ with $|d_i(\mathbb{G}(W,n))/(n-1)-d|\geq\eps$ is at
most $2ne^{-(n-1)\eps^2/2}$.
\end{proof}

\begin{definition}
We say that $\{G_i\}_{i=1}^\infty$ is a {\it degree-uniformly}
convergent sequence of simple graphs with limiting degree $0\leq
d\leq 1$ if it is convergent, $\lim_{i\to\infty}|V(G_i)|=\infty$ and
\[
\lim_{i\to\infty}\frac{d_{\rm max}(G_i)}{|V(G_i)|} =
\lim_{i\to\infty}\frac{d_{\rm min}(G_i)}{|V(G_i)|}=d.
\]
\end{definition}

\begin{lemma}\label{reglim}
If a CR-graphon $W$ is $d$-regular then there is a sequence of
CR-graphs $\{G_n\}_{n=1}^\infty$ that degree-uniformly converges to
$W$.
\end{lemma}

\begin{proof}
The Borel-Cantelli lemma together with lemma \ref{degcont} implies
that with probability one
\[
\lim_{i\to\infty} \frac{d_{\rm max}(\mathbb{G}(W,n))}{n-1}
=\lim_{i\to\infty} \frac{d_{\rm min}(\mathbb{G}(W,n))}{n-1}=d.
\]
It is clear that with probability one $\mathbb{G}(W,n)$ does not
contain an induced $P_4$, and hence it is a CR-graph for every $n$.
We also know \cite{LSz1} that with probability $1$ it converges to
$W$ as $n\to\infty$. These facts imply that the sequence
$\{\mathbb{G}(W,n)\}_{n=1}^\infty$ satisfies the conditions with
probability $1$.
\end{proof}

\begin{lemma}\label{felsobecskomp}
Let $G$ be a finite disconnected CR-graph on $n\ge 2$ nodes such that
$|d(v)/n-d|\leq\eps$ for every $v\in V(G)$ with some $d\in[0,1]$.
Then every connected component of $G$ has size less than
$(\frac23+\frac43\eps)n$.
\end{lemma}

\begin{proof}
Let $G'$ be a connected component of $G$ of maximal size. Let
$n=|V(G)|$ and $a=|V(G')|$. We may assume that $a>1$ (else, the
assertion is trivial). Let $v\in V(G)\setminus V(G')$, then $d(v)
<n-a$. Since $G'$ is a connected CR-graph, there is a node $w\in
V(G')$ with degree at least $a/2$. Then by our assumption
\[
2\eps n\geq d(w)-d(v) > \frac{a}2-(n-a)= \frac32 a-n.
\]
This implies that $a<(\frac23+\frac43\eps)n$.
\end{proof}

\begin{lemma}\label{LEM:DIRSUM}
Let $W$ be a $d$-regular CR-graphon. Then either $W$ or $1-W$ can be
decomposed as a weighted direct sum of at least two regular
CR-graphons.
\end{lemma}

\begin{proof}
By Lemma \ref{reglim}, there is a sequence of CR-graphs
$\{G_n\}_{n=1}^\infty$ that degree-uniformly converges to $W$. For
each $n$, either $G_n$ or $\overline{G}_n$ is disconnected, since
$G_n$ is a CR-graph. We may assume, by restricting ourselves to a
subsequence, that either $G_n$ is disconnected for all $n$, or
$\overline{G}_n$ is disconnected for all $n$. By complementing if
necessary, we may assume that $G_n$ is disconnected for all $n$.

Let $H_{n,1},\dots,H_{n,k_n}$ be the connected components of $G_n$.
Since the convergence is degree-uniform, it follows that for any
$0<d'<d$, all degrees of $G_n$ are larger than $d'|V(G_n)|$ if $n$ is
large enough, and then trivially $|V(H_{n,i})|\ge d'|V(G_n)|$. This
implies that $k_n$ remains bounded, and so by going to a subsequence
again, we may assume that $k_n=k$ is independent of $n$. By the same
token, we may assume that $|V(H_{n,i})|/|V(G_n)|$ has a limit $a_i$
as $n\to\infty$. Clearly $a_i\ge d$ and $\sum_i a_i=1$. We may also
assume that for each $1\le i\le k$, the sequence of graphs
$(H_{n,i})_{n=1}^\infty$ is convergent. Let $W_i$ denote its limit
graphon. It is straightforward to check that $W_i$ is a regular
CR-graphon. Furthermore, the weighted direct sum $\bigoplus_i (a_i)
W_i$ is the limit of the graphs $G_n$. By the uniqueness of the
limit, $W$ is weakly isomorphic to $\bigoplus_i (a_i) W_i$.
\end{proof}

Now we are able to prove Theorem \ref{THM:REG-LOCFIN}.

\begin{proof}
Assume that $W=\bigoplus_{i=1}^k (a_i) W_i$, $k\ge 2$. We may assume
that the $W_i$ cannot be written as weighted direct sums in a
nontrivial way. We build a tree by starting with a root corresponding
to $W$, having $k$ children corresponding to $1-W_1,\dots,1-W_k$. If
any of these functions is almost everywhere $0$, then this node will
be a leaf. Else, we continue building the tree from this node as
root.

If $W$ cannot be written as a weighted direct sum of at least two
regular CR-graphons, then by Lemma \ref{LEM:DIRSUM}, $1-W$ can be,
and we start the tree with a root with a single child, corresponding
to $1-W$.

This way we obtain a tree $T$, where each node $v$ is labeled by a
regular CR-graphon $W_v$. For each node of the tree constructed this
way, we define $f(v)$ as the product of the weights of the graphons
along the path from the root to $v$. It is straightforward to check
that $W$ is weakly isomorphic to the graphon $U_T$ with the
probability distribution defined by $f$.
\end{proof}

Let $W$ be a regular CR-graphon represented by the tree $T$ with a
measure $\mu$ on $\Omega_T$.
For a node $v\in V(T)$ let
$f(v):=\mu(\Omega_v)$. We assume that $f(v)>0$, since parts of the tree
with $0$ weight can be deleted.
We observe that for every $v\in V(T)$, the
subtree $T_v$ of $T$ rooted at $v$, with the same local distributions
(i.e., with node weights $f(u)/f(v)$), defines another regular
CR-graphon. This implies that the value
\[
c(v)=\int\limits_{\Omega_v} U_{T_v}(x,y)\,d\mu(y)
\]
is the same for all $x\in\Omega_v$. The degree of the graphon on
$T_v$ is $d(v)=c(v)/f(v)$. (Note however that, depending on the
parity of the depth of $v$, either $U_{T_v}$ or $1-U_{T_v}$ is an
induced subgraphon of $U_T$.) For $u\in V(T)$ and $v\in C_u$ we have
\begin{equation}\label{EQ:CCA}
c(u)+c(v)=f(v),
\end{equation}
and for every leaf $u$ (if any)
\begin{equation}\label{EQ:CCA2}
c(u)=0.
\end{equation}

The following simple lemma gives some conditions that $f$ and $c$
satisfy.

\begin{lemma}\label{LEM:BOUNDS}
Let $W_T$ be a regular CR-graphon.

(a) If $u\in V(T)$ has $r$ children, then $c(u)\le \frac1r f(u)$.

(b) If $u\in V(T)$, $v\in C_u$ and $v$ has $r$ children, then
$f(u)\ge (2-\frac1r) f(v)$.
\end{lemma}

\begin{proof}
(a) Let $v_1,v_2\dots,v_r$ be the children of $u$. By \eqref{EQ:CCA},
$f(v_i)=c(u)+c(v_i) \ge c(u)$, and summing this over $i$, we get
$f(u) \ge  r c(u)$.

(b) Let $v'$ be a sibling of $v$. Using \eqref{EQ:CCA} and (a),
\[
f(v')=c(u)+c(v') \ge c(u) = f(v)-c(v) \ge (1-\frac1r)f(v),
\]
and so
\[
f(u)\ge f(v)+f(v') \ge(2-\frac1r)f(v).
\]
\end{proof}

\begin{lemma}\label{LEM:CF}
Let $T$ be a locally finite tree such that no node except possibly
the root has exactly one child. Let $c,f:~V(T)\to\R_+$ be two
functions satisfying \eqref{EQ:F-ADD}, \eqref{EQ:CCA} and
\eqref{EQ:CCA2}. Then the probability measure defined by $f$ gives a
regular CR-graphon on $T$.
\end{lemma}

\begin{proof}
By \eqref{EQ:F-ADD}, the function $f$ defines a probability measure
$\pi$ on $(\Omega,\AA)$. Let $x=(v_0,v_1,v_2,\dots)$ be a maximal path
starting at the root $r=v_0$. Lemma \ref{LEM:BOUNDS}(b) implies that
if this path is infinite, then $f(v_n)\to 0$, and part (a) of the
same Lemma implies that $c(v_n)\to 0$.

The path $x$ is connected to all paths $y$ that branch off from $x$
at $v_1$, $v_3$, $v_5$, etc. The $\pi$-measure of these paths is
$(f(v_1)-f(v_2))+(f(v_3)-f(v_4))+\dots$, which by \eqref{EQ:CCA} can
be written as
\[
(c(v_0)+c(v_1))-(c(v_1)+c(v_2))+(c(v_2)+c(v_3))-\dots = c(v_0)
\]
(if the path ends at a leaf, then we use \eqref{EQ:CCA2}). This is
indeed independent of the path.
\end{proof}

Now we are able to prove the second main result in this section:

\begin{theorem}\label{LEM:LOCFIN-CRREG}
For every locally finite rooted tree $T$ there is a unique regular
CR-graphon on $T$.
\end{theorem}

\begin{proof}
{\it Existence.} First we prove this for a finite tree, by induction
on the depth. For a single node, the function $U\equiv0$ is a regular
CR-graphon.

Suppose that the tree has more than one node, and let $u_1,\dots,u_k$
be the children of the root. By induction, we find regular
CR-graphons on $T_{u_1},\dots,T_{u_k}$, with degrees $d_1,\dots,d_k$.
Note that since $u_i$ is either a leaf or has at least two children,
we must have $d_i<1$. Let
\[
d=\frac1{\sum_i \frac1{1-d_i}},   \qquad a_i=\frac{d}{1-d_i},
\]
then scaling the measure of $T_{u_i}$ by $a_i$, complementing each
$T_{u_i}$ and taking their disjoint union, we get a $d$-regular
CR-graphon on $T$.

Now suppose that $T$ is infinite, and let $T_k$ denote the tree
obtained by deleting all nodes farther than $k$ from the root. By the
above, there is a regular CR-graphon on $T_k$, which yields two
functions $f^k$ and $c^k$ on $V(T_k)$ satisfying \eqref{EQ:F-ADD},
\eqref{EQ:CCA} and \eqref{EQ:CCA2}. We can select a subsequence of
the indices $k$ such that $f^k(v)$ tends to some $f(v)$ and $c^k(v)$
tends to some $c(v)$ as $k$ ranges through this subsequence. Clearly,
the functions $f$ and $c$ also satisfy \eqref{EQ:F-ADD},
\eqref{EQ:CCA} and \eqref{EQ:CCA2}, and so by Lemma \ref{LEM:CF},
they yield a regular CR-graphon on $T$.

\smallskip

{\it Uniqueness.} Similarly to the existence part, it is easy to see
that uniquness holds for finite trees. If the tree is infinite then
it contains an infinite path $P=(v_0,v_1,v_2,\dots)$ starting at the
root. Suppose that there are two weightings $f,f'$ of the nodes of
$T$ such that they both define a regular CR-graphon, then we get by
\eqref{EQ:CCA}
\[
f(v_1)-f(v_2)+\dots+(-1)^{k+1}f(v_k) = c(v_0)+(-1)^{k+1}c(v_k).
\]
The sequence $f(v_k)$ is monotone decreasing and it tends to $0$.
Hence $c(v_k)\to 0$, and
\[
c(v_0)=\sum_{k=1}^\infty (-1)^{k+1}f(v_k).
\]
In the other weighting,
\[
c'(v_0)=\sum_{k=1}^\infty (-1)^{k+1}f'(v_k).
\]
Let
\[
z_i=\frac{f'(v_{i})f(v_{i-1})}{f'(v_{i-1})f(v_i)},
\]
then $f'(v_k)= z_1\dots z_k f(v_k)$. Thus
\begin{equation}\label{EQ:CPRIME}
c'(v_0)=\sum_{k=1}^\infty (-1)^{k+1} z_1\cdots z_kf(v_k) =
z_1\bigl(f(v_1)-z_2(f(v_2) -\dots)\bigr).
\end{equation}

Choose the path $v_0,v_1,\dots$ as follows. Given $v_i$, choose
$v_{i+1}\in C_{v_i}$ so that $z_{i+1}$ is as large as possible if $i$
is odd, and as small as possible if $i$ is even. Clearly
$z_{i+1}\ge1$ if $i$ is odd, and $z_{i+1}\le 1$ if $i$ is even. This
is clearly possible. Raising $z_1$ to $1$, then lowering $z_2$ to 1,
then raising $z_3$ to $1$ etc. increases the expression in
\eqref{EQ:CPRIME}, and hence
\begin{equation}\label{EQ:CFFC}
c'(v_0)\le c(v_0).
\end{equation}
Since the reverse inequality follows similarly, we get that
$c'(v_0)=c(v_0)$. From the fact that equality holds in
\eqref{EQ:CFFC}, we get that all $z_1=1$, which in turn implies that
$f'(v_1)=f(v_1)$ for any child $v_1$ of $v_0$.

Applying the same argument to the graphons $W_{T_{v_1}}$, $v\in
C_{v_0}$, then to their children etc., we get that $f(v)=f'(v)$ for
all $v$.
\end{proof}

\subsection{Forcible regular CR-graphons with irrational edge densities}
\label{SEC:IRR}

We start with an easy observation:

\begin{lemma}\label{LEM:STEP-REG}
If a $d$-regular CR-graphon is a stepfunction, then $d$ is rational.
\end{lemma}

\begin{proof}
By induction on the depth of the tree.
\end{proof}

In contrast to this, we prove that for every $\alpha\in[0,1]$ there
is a regular CR-graphon of degree $\alpha$ which is finitely
forcible. As lemma \ref{LEM:STEP-REG} shows, for irrational $\alpha$
such a graphon in not a step function.

Let $P_3^2$ be the disjoint union of two copies of $P_3$. Let
$\mathcal{Z}$ denote the set of regular CR-graphons that don't
contain any induced copy of $P_3^2$ and its complement.

\begin{lemma}
The set $\mathcal{Z}$ consists of those regular CR-graphons whose
representing tree has the property that every node has at most one
child which is not a leaf.
\end{lemma}

\begin{proof}
Let $W$ be an element in $\mathcal{Z}$. First of all note that a
graphon without an induced copy of $P_3$ is the disjoint union of
complete graphons. It follows that a graphon without an induced copy
of $P_3^2$ has at most one connected component which is not complete.
Applying this for the sub-graphons (and their complements)
corresponding to the nodes of the tree of $W$ we get that $W$
satisfies the condition. The other direction is trivial.
\end{proof}

The previous lemma shows that the tree of an element in $\mathcal{Z}$
is one (possibly infinite) path with additional leaves hanging from
its nodes. Thus the structure is determined by the integer sequence
$n_1,n_2,\dots$ where $n_k$ is the number of leaves at the $k$-th
level.

We start with a simple example. Let $\alpha=(3-\sqrt{5})/2$. Note
that $\alpha=(1-\alpha)^2$. There exists a unique graphon $W$ which
is the disjoint union of a clique of size $\alpha$ and a version of
the complement of $W$ scaled to the size $1-\alpha$. The graphon $W$
has an iterated structure. The choice of $\alpha$ guarantees that $W$
is $\alpha$-regular. We show that $W$ is finitely forcible.

\begin{lemma}
The graphon $W$ is the only element of $\mathcal{Z}$ with degree
$\alpha$ and thus it is finitely forcible.
\end{lemma}

\begin{proof}
Let $W'$ be another graphon with the above property. Since $1/3
<\alpha < 1/2$, we have that $W'=(\alpha)1 \oplus (1-\alpha)W''$,
where $W''$ has density $\alpha/(1-\alpha)=1-\alpha$.
This shows that we can
inductively continue the process with the complement $W''$ and
finally obtain the desired iterated structure.
\end{proof}

The previous argument can be easily generalized to other irrational
values of $\alpha$.

\begin{prop}\label{itstr1}
For every irrational number $0<\alpha<1$ there is exactly one graphon
in $\mathcal{Z}$ with edge density $\alpha$ and so this graphon is
finitely forcible.
\end{prop}

\begin{proof}
Let $\alpha$ be a number between $0$ and $1/2$ (the case $\alpha
>1/2$ is similar). Let $W$ be a graphon from $\mathcal{Z}$ with edge
density $\alpha$. Since $\alpha<1/2$, the graphon $W$ is weighted
direct sum of $n_1$ cliques, all with weight $\alpha$ and a connected
element $W'$ of $\mathcal{Z}$ with weight $0<s<1$ and with edge
density between $1/2$ and $1$. To guarantee edge density $\alpha$ in
this component, $s$ has to be between $\alpha$ and $2\alpha$.
Consequently $n_1$ is the unique natural number such that $\alpha
\leq 1-n_1\alpha<2\alpha$. The complement $1-W'$ is an element from
$\mathcal{Z}$ with edge density smaller than $1/2$ and we can
continue the process.
\end{proof}

The reader can see that the number $n_1,n_2,\dots$ are uniquely
determined by $\alpha$ so it is a natural question to ask what these
numbers are. An elementary calculation shows that these numbers are
basically the numbers occurring in the (unique) continued fraction
expansion of $\alpha$. The only exception is the first number, which
is shifted by one:
\[
\alpha=\cfrac{1}{n_1+1
+\cfrac{1}{n_2+\cfrac{1}{n_3+\cfrac{1}{\ddots}}}}.
\]
This shows that one can force graphons that encode an arbitrary
sequence of natural numbers in a very structural way.

\subsection{Forcing the binary tree}

In this section we prove that the regular CR-graphon $U_{T_2}$  is
finitely forcible, where $T_2$ is the complete binary tree of
infinite depth, with a root of degree $1$ added to comply with our
previous definitions. We note that $U_{T_2}$ has the following
alternative description: Consider the space $V(C_4)^{\N}$ with the
uniform probability measure. We connect two nodes $x$ and $y$ of
$V(C_4)^{\N}$ if for the first coordinate where they differ, say
$i\in\N$, $x_i$ and $y_i$ are connected in $C_4$. The graphon
$U_{T_2}$ can be called the {\it infinite lexicographic power} of
$C_4$.

Let us define the following signed labeled graphs: $B$ is $K_3$ with
one node labeled $1$, the incident edges signed ``$+$'', and the
opposite edge signed ``$-$''; and $C$ and $D$ are obtained from
$K_2^\bullet$ and $B$, respectively, by adding a new node labeled
$2$, and connecting it to all other nodes by edges signed
``$-$''. Also consider the signed graphs $\overline{B}$,
$\overline{C}$ and $\overline{D}$ obtained from $B$, $C$ and $D$ by
switching the ``$+$'' and ``$-$'' signs on the edges.

\begin{prop}\label{PROP:TREE-FORC}
Let $W$ be a graphon satisfying
\[
t(\widehat{P}_4,W)=0,\qquad t_1(K_2^\bullet,W)(x)=\frac{2}{3},
\]
\[
t_1(B,W)(x)=\frac{8}{45},~~t_1(\overline{B},W)(x)=\frac{2}{45}
\]
and
\[
2t_2(C,W)^2(x,y)=5t_2(D,W)(x,y),
~~2t_2(\overline{C},W)^2(x,y)=5t_2(\overline{D},W)(x,y)
\]
almost everywhere. Then $W$ is weakly isomorphic to $U_{T_2}$.
\end{prop}

\begin{proof}
It is straightforward to check that the graphon $U_{T_2}$ satisfies
these identities.

The first two identities mean that $W$ is a regular CR-graphon with
degree $2/3$. By Theorem \ref{THM:REG-LOCFIN} we know that $W$ can be
represented by a locally finite tree $T$ and so we can assume that
$W=U_T$. The edge density $2/3$ guarantees that the root $r$ of $T$
has one child $q$, but $q$ must have at least $2$ children, i.e.,
$U_T$ is a connected graphon and $1-U_T$ has at least $2$ components.
Let $v$ be any child of $q$, and let $\Omega'=\Omega\setminus
\Omega_v$. We also know that $v$ has at least two children $v_1,v_2$.
Let $x\in\Omega_{v_1}$ and $y\in\Omega_{v_2}$. By the definition of
$U_T$, $U_T(x,z)=0$ for $z\in \Omega_{v_2}$, $U_T(y,z)=0$ for $z\in
\Omega_{v_1}$and $U_T(x,z)=U_T(y,z)=1$ for $z\in \Omega'$. Hence
\begin{align}\label{TXYCUT}
t_2(C,U_T)(x,y)&=\int\limits_{\Omega} U_T(x,z)(1-U_T(z,y))\,dz =
\int\limits_{\Omega_{v_1}} U_T(x,z)\,dz
\nonumber\\
&= \int\limits_{\Omega} U_T(x,z)\,dz - \mu(\Omega') =
t_1(K_2^\bullet,U_T)(x)- \mu(\Omega') = \frac23 - \mu(\Omega').
\end{align}
Similarly,
\begin{align*}
t_2(D,U_T)(x,y)&=\int\limits_{\Omega\times\Omega}
U_T(x,z)(1-U_T(z,y))U_T(x,u)(1-U_T(u,y))(1-U_T(z,u))\,du\,dz\\
&= \int\limits_{\Omega_{v_1}\times\Omega_{v_1}}
U_T(x,z)U_T(x,u)(1-U_T(z,u))\,du\,dz\\
&= t_1(B,U_T)(x)-\int\limits_{\Omega'\times\Omega'}
(1-U_T(z,u))\,du\,dz.
\end{align*}
Here
\[
\int\limits_{\Omega'\times\Omega'} 1-U_T(z,u)\,du\,dz =
\int\limits_{\Omega\times\Omega'} 1-U_T(z,u)\,du\,dz =
\int\limits_{\Omega'}
1-t_1(K_2^\bullet,U_T)(u)\,du=\frac13\mu(\Omega'),
\]
and so
\begin{equation}\label{TXYDUT}
t_2(D,U_T)(x,y)=\frac{8}{45}-\frac13\mu(\Omega').
\end{equation}
Using our conditions, we get from \eqref{TXYCUT} and \eqref{TXYDUT}
\[
2\Bigl(\frac23-\mu(\Omega')\Bigr)^2=2t_2(C,U_T)(x,y)^2
=5t_2(D,U_T)(x,y)=\frac89-\frac53\mu(\Omega').
\]
This simplifies to the equation $2\mu(\Omega')^2=\mu(\Omega')$. Since
$\mu(\Omega')\not=0$ (as $q$ has at least two children), we get
$\mu(\Omega')=1/2$ and so $\mu(\Omega_v)=1/2$. This is true for every
child of $q$, and hence there are exactly two children, both with
weight $1/2$.

To finish, it is easy to check that the complement of the graphon
$U_{T_v}$ satisfies the same identities as listed in the statement
for either child $v$ of $q$. Iterating the argument, we get that $T$
is a complete binary tree.
\end{proof}

\begin{corollary}\label{UB-FORCE}
The graphon $U_{T_2}$ is finitely forcible.
\end{corollary}

\begin{proof}
Lemma \ref{LEM:UNLABEL2} guarantees that in Proposition
\ref{PROP:TREE-FORC} the first four equations are forcible by a
finite number of subgraph densities, but for the last two equations
we need a simple additional argument. Let $G=(2C^2-5D)^2$ and
$H=(2\overline{C}^2-5\overline{D})^2$ in the algebra of $2$-labeled
quantum graphs (with multiple edges allowed). It is clear that
$t(\lunl G\runl ,W)=0$ and $t(\lunl H\runl ,W)=0$ force the desired
equations in $\WW$. Let $G_2$ and $H_2$ be the quantum graphs
obtained from $\lunl G\runl $ and $\lunl H\runl $ by reducing
multiple edges. As noted above, the first two conditions in
Proposition \ref{PROP:TREE-FORC} imply that $W$ is regular, and by
Theorem \ref{THM:REG-LOCFIN}, $W$ is $0$-$1$ valued. Within the set
of $0$-$1$ valued graphons this edge reduction doesn't change the
density functions, so $t(G_2,W)=0$ and $t(H_2,W)=0$ are equivalent to
$t(\lunl G\runl ,W)=0$ and $t(\lunl H\runl ,W)=0$.
\end{proof}

\section{Necessary conditions for finite forcing}

\subsection{Infinite rank}

We define the {\it rank} of a graphon $W$ as its rank as a kernel
operator. In other words, the rank of $W$ is the least nonnegative
integer $r$ such that there are measurable functions
$w_i:~[0,1]\to\R$ and reals $\lambda_i$ ($i=1,\dots,r$) such that
\begin{equation}\label{WRANK}
W(x,y)=\sum_{k=1}^r \lambda_k w_k(x)w_k(y)
\end{equation}
almost everywhere. If no such integer $r$ exists, then we say that
$W$ has infinite rank.

\begin{theorem}\label{FINRANK}
If $W$ has finite rank, then for every finite list $F_1,\dots,F_m$ of
simple graphs there is a stepfunction $U$ such that
$t(F_i,U)=t(F_i,W)$ for $i=1,\dots,m$.
\end{theorem}

\begin{proof}
We know that $W$ has a decomposition \eqref{WRANK}, where we may
assume that the $\lambda_k$ are the eigenvalues and the $w_k$ are the
corresponding eigenfunctions of $W$ as a kernel operator. In this
case, all the $w_k$ are bounded and all the moments
$M(\{w_1,\dots,w_r\},k)$ are finite.

Fix a simple graph $F=(V,E)$. For a map $\varphi:~E\to[r]$,
$t\in[r]$, and $i\in V$, let $d_t(\varphi,i)$ denote the number of
edges $e\in E$ incident with $i$ for which $\varphi(e)=t$, and set
$\lambda_\varphi=\prod_{ij\in E}\lambda_{\varphi(ij)}$. Then
\begin{align*}
t(F,W)&=\int\limits_{[0,1]^V} \prod_{ij\in E} W(x_i,x_j)\,dx
=\int\limits_{[0,1]^V} \prod_{ij\in E}\left(\sum_{k=1}^r \lambda_k
w_k(x_i)w_k(x_j)\right)\,dx \\
&=\int\limits_{[0,1]^V} \sum_{\varphi\in [r]^E}
\lambda_\varphi\prod_{ij\in
E} w_{\varphi(ij)}(x_i)w_{\varphi(ij)}(x_j)\,dx\\
&=\int\limits_{[0,1]^V} \sum_{\varphi\in [r]^E} \lambda_\varphi
\prod_{i\in V} \prod_{t\in[r]} w_t(x_i)^{d_t(\varphi,i)}\,dx
=\sum_{\varphi\in [r]^E} \lambda_\varphi \prod_{i\in V}
\int\limits_0^1 \prod_{t\in[r]} w_t(y)^{d_t(\varphi,i)}\,dy
\\&=\sum_{\varphi\in [r]^E} \lambda_\varphi \prod_{i\in V}
M(w,d(\varphi,i)).
\end{align*}
So if $(u_1,\dots,u_r)$ is another set of functions that satisfy
\begin{equation}\label{UW-MOM}
M(u,d(\varphi,i))=M(w,d(\varphi,i))
\end{equation}
for every $1\le j\le m$, $i\in V(F_j)$ and $\varphi:~V(F_j)\to[r]$,
then the function
\[
U=\sum_{t=1}^r \lambda_t u_t(x) u_t(y)
\]
satisfies $t(F_j,U)=t(F_j,W)$ for all $j=1,\dots,m$. By Theorem
\ref{ONE-MOM} there is a system of functions $u$ satisfying
\eqref{UW-MOM} which are stepfunctions, and then $U$ is also a
stepfunction.
\end{proof}

\begin{corollary}\label{INF-RANK}
Every finitely forcible graphon is either a stepfunction or has
infinite rank.
\end{corollary}

In view of Theorem \ref{MON-POLY}, the following corollary of this
theorem may be surprising:

\begin{corollary}\label{POLY-NOT}
Assume that $W\in\WW_0$ can be expressed as a non-constant polynomial
in $x$ and $y$. Then $W$ is not finitely forcible.
\end{corollary}

\subsection{Weak homogeneity}

For every graph $F=(V,E)$, and node $i\in V$, let $F^i$ denote the
$1$-labeled quantum graph obtained by labeling $i$ by $1$, and for
every edge $ij\in E$, let $F^{ij}$ denote the 2-labeled quantum graph
obtained from $F$ by deleting the edge $ij$, and labeling $i$ by $1$
and $j$ by $2$. Let $F^\dag=\sum_{i\in V} F^i$ and
$F^{\ddag}=\frac{1}{2}\sum_{i,j:\,ij\in E} F^{ij}$ (each edge
contributes two terms, since its endpoints can be labeled in two
ways). We extend the operators $F\to F^\dag$ and $F\to F^{\ddag}$
linearly to all quantum graphs.

\begin{example}\label{EXA:CIRC}
Clearly $C_n^{\ddag}=nP_n^{\bullet\bullet}$, where
$P_n^{\bullet\bullet}$ denotes the path on $n$ nodes with its
endpoints labeled. So
\[
t_2(C_n^\ddag,W) = n W^{\circ(n-1)}.
\]
\end{example}

These operations were introduced by Razborov \cite{Razb1,Razb2} in
the proof of conjectures about the minimum number of triangles in
simple graphs with given edge density. For us, their significance is
in the following formulas.

We consider $\WW$ as a Banach space with the $L_\infty$ norm. Let
$U_t$, $0\le t\le 1$ be a family of functions in $\WW$. We say that
$U_t$ is {\it differentiable} if for every $t\in[0,1]$ there exists a
function $\dot{U}_t\in\WW$ such that
\[
\Bigl\|\frac1{s-t}(U_s-U_t)-\dot{U}_t\Bigr\|_\infty \to 0 \qquad
(s\in[0,1],~s\to t).
\]

\begin{lemma}\label{LEM:WDIFF}
Let $U_t$, $0\le t\le 1$ be a uniformly bounded differentiable family
of functions in $\WW$ and $F=(V,E)$, a simple graph. Then the
function $t(F,U_t)$ is differentiable as a function of $t$, and
\[
\frac{d}{dt} t(F,U_t) = \bigl\langle \dot{U}_t, t_2(F^{\ddag},U_t)
\bigr\rangle.
\]
\end{lemma}

\begin{proof}
Suppose that $\|U_s\|_\infty\le C$ for some real number $C$. Write
\begin{equation}\label{EQ:TUH}
t(F,U_{t+h})-t(F,U_t)= \bigl(t(F,U_t+h\dot{U}_t)-t(F,U_t)\bigr)+
\bigl(t(F,U_{t+h})-t(F,U_t+h\dot{U}_t)\bigr).
\end{equation}
Here the first term is a polynomial in $h$:
\[
t(F,U_t+h\dot U_t)-t(F,U_t) = h\langle\dot{U}_t,
t_2(F^\ddag,U_t)\rangle + O(h^2),
\]
while the second term can be estimated by \eqref{EQ:COUNT} and the
definition of differentiation:
\[
|t(F,U_{t+h})-t(F,U_t+h\dot{U}_t)|\le
|E(F)|C^{|E(F)|-1}\|U_{t+h}-U_t-h\dot{U}_t\|_\infty = o(h).
\]
So indeed
\[
\frac1h(t(F,U_{t+h})-t(F,U_t))\to \langle\dot{U}_t,
t_2(F^\ddag,U_t)\rangle \qquad (h\to 0).
\]
\end{proof}

We use this formula to derive a necessary condition for finite
forcibility.

\begin{lemma}\label{LEM:NONFIN-DIFF}
Suppose that $W\in\WW$ is forced (in $\WW$) by the simple graphs
$F_1,\dots,F_m$. Also suppose that there exists an open ball
$\UU\subseteq\WW$ about $W$ and a Lipschitz map $\Phi:~\UU\to \WW$
such that for all $U\in\UU$ the function $\Phi(U)$ satisfies $\langle
\Phi(U), t_2(F_i^\ddag,U) \rangle=0$ $(i=1,\dots,m)$. Then $\langle
\Phi(W), t_2(F^\ddag,W) \rangle=0$ for every simple graph $F$.
\end{lemma}

\begin{proof}
By classical results on differential equations in Banach spaces (see
e.g.~\cite{Zeid}), there exists a $b>0$ and a differentiable family
$\{U_s:~s\in[-b,b]\}$ of functions in $\UU$ satisfying the
differential equation
\[
\dot{U}_s=\Phi(U_s), \quad U_0=W.
\]

Lemma \ref{LEM:WDIFF} shows that for every simple graph $F$
\[
\frac{d}{ds}t(F,U_s)=\Bigl\langle\dot{U}_s,
t_2(F^\ddag,U_s)\Bigr\rangle= \langle\Phi(U_s),
t_2(F^\ddag,U_s)\rangle.
\]
In particular, we have
\[
\frac{d}{ds}t(F_i,U_s)= \langle\Phi(U_s), t_2(F_i^\ddag,U_s)\rangle=0
\]
for $i=1,\dots,m$, and hence $t(F_i,U_s)=t(F_i,U_0)=t(F_i,W)$ for all
$\in[0,c]$. Since the graphs $F_i$ force $W$, it follows that the
$U_s$ are weakly isomorphic to $W$, and so $t(F,U_s)=t(F,W)$ for
every $F$. But then $\langle \Phi(W), t_2(F^\ddag,W)
\rangle=\frac{d}{ds}t(F,U_s)\bigl|_{s=0}=0$ as claimed.
\end{proof}

Let $\LL(W)$ be the linear space generated by 2-variable functions
$t_2(F^\ddag,W)$ (modulo zero sets). Inequality (\ref{EQ:COUNT2})
implies that $t_2(F^\ddag,W)\in \WW$ for all $F$. Due to the identity
\begin{equation}\label{EQ:F1F1DDAG}
t_2((F_1F_2)^\ddag,W)=  t(F_1,W) t_2(F_2^\ddag,W) + t(F_2,W)
t_2(F_1^\ddag,W),
\end{equation}
the space $\LL(W)$ is generated by functions $t_2(F^\ddag,W)$ where
$F$ is connected.

\begin{prop}\label{PROP:LINTF}
Let $W$ be a finitely forcible graphon. Then $\LL(W)$ has finite
dimension if and only if $W$ is a stepfunction.
\end{prop}

\begin{proof}
If $W$ is a stepfunction, then every function $t_2(F^\ddag,W)$ is a
stepfunction with the same steps, and so $\LL(W)$ is finite
dimensional. Conversely, if $\LL(W)$ is finite dimensional, then by
Example \ref{EXA:CIRC}, the functions $W^{\circ k}\in\LL(W)$ are
linearly dependent, and so $W$ satisfies a polynomial equation as an
operator. This means that it has a finite number of different nonzero
eigenvalues. Since every nonzero eigenvalue has finite multiplicity,
$W$ has finite rank. By Corollary \ref{INF-RANK}, $W$ is a
stepfunction.
\end{proof}

\begin{remark}\label{INFINITESIMAL}
Suppose that $\LL(W)$ has finite dimension and the functions
$t_2(F_1^\ddag,W),\dots,t_2(F_k^\ddag,W)$ generate it. Informally,
this means that every infinitesimal change in $W$ that preserves
$t(F_1,W),\dots,t(F_k,W)$, also preserves $t(F,W)$ for every $F$; we
could say that $W$ is {\it infinitesimally finitely forcible}.
Proposition \ref{PROP:LINTF} says that graphons that are both
finitely forcible and infinitesimally finitely forcible are exactly
the stepfunctions.

Our examples of finitely forcible non-step-functions (e.g.,
half-graphs) show that there are graphons which are finitely forcible
but not infinitesimally finitely forcible. We don't know if the
converse is true.
\end{remark}

\begin{lemma}\label{LEM:FIN-NEC}
Suppose that $W\in\WW$ is forced (in $\WW$) by the simple graphs
$F_1,\dots,F_m$. Then either $t_2(F_1^\ddag,W),\dots,
t_2(F_m^\ddag,W)$ are linearly dependent, or they generate $\LL(W)$.
\end{lemma}

\begin{proof}
Suppose not, then there is a simple graph $F_{m+1}$ such that
$t_2(F_1^\ddag,W),\dots, t_2(F_m^\ddag,W)$ and $t_2(F_{m+1}^\ddag,W)$
are linearly independent. For $U\in\WW$, set
$h_k(U)=t_2(F_k^\ddag,U)$. Let $\Phi(U)$ denote the component of
$h_{m+1}(U)$ orthogonal to the subspace spanned by
$h_1(U),\dots,h_m(U)$.

We need the following technical claim.

\begin{claim}\label{CLAIM:LIP}
There is an open ball $\UU$ in $\WW$ centered at the function $W$
such that $\Phi:~\WW\to\WW$ is Lipschitz on $\UU$.
\end{claim}

Inequality \eqref{EQ:COUNT2} implies that there is a neighborhood
$\UU$ of $W$ in $\WW$ such that the functions
$h_1(U),\dots,h_{m+1}(U)$ are linearly independent for $U\in\UU$. We
may assume that $U$ is an open ball such that $\|U\|_\infty\le
2\|W\|_\infty$ for all $U\in\UU$.

Let $(g_1(U),\dots,g_{m+1}(U))$ be the Gram-Schmidt orthogonalization
of $(h_1(U),\dots,h_{m+1}(U))$, then $\Phi(U)=g_{m+1}(U)$. For a
function $H\in\WW$, let $\Psi_H:~\WW\to\WW$ denote the orthogonal
projection onto the subspace orthogonal to $H$. Consider the
functions
\[
g_{k,r}(U)=\Psi_{g_r(U)}\dots\Psi_{g_1(U)}h_k(U).
\]
Then we have
\[
g_{k,r+1}(U)=\Psi_{g_{r+1}(U)} g_{k,r}(U),
\]
and
\[
g_{k,0}(U)=h_k(U), \qquad g_{k+1,k}(U)=g_{k+1}(U).
\]

We prove by induction on $k$ and $r$ ($r<k$) that there is a constant
$c_{k,r}>0$ and an open ball $\UU_{k,r}$ about $W$ such that if
$\{U_s:~0\le s\le 1\}\subseteq\UU_{k,r}$ is a differentiable family,
then
\[
\Bigl\|\frac{d}{ds} g_{k,r}(U_s)\Bigr\|_\infty\le c_{k,r}
\Bigl\|\frac{d}{ds} U_s\Bigr\|_\infty.
\]
This will imply that $\Phi = g_{m+1}$ is Lipschitz.

First, for the functions $g_{k,0}(U)=h_k(U)=t_2(F_k^\ddag,U)$ this
follows from inequality \eqref{EQ:COUNT2}.

Let $k\ge 1$, $r<k$, and suppose that we know the existence of
$c_{r+1,0}$ and of $c_{k+1,r}$. We prove that $c_{k+1,r+1}$ exists.

Set $G=g_{k+1,r}(U_s)$ and $H=g_{r+1}(U_s)$, then (denoting
differentiation by dot)
\begin{align}\label{DPSI}
\frac{d}{ds} \Psi_H G &= \frac{d}{ds} \Bigl( G - \frac{\langle G,H
\rangle}{\langle H,H \rangle}H\Bigr)\nonumber\\
&= \dot{G} - \frac{\langle \dot{G},H \rangle}{\langle H,H \rangle}H -
\frac{\langle G,\dot{H} \rangle}{\langle H,H \rangle}H + 2
\frac{\langle G,H \rangle\cdot\langle H,\dot{H} \rangle}{\langle H,H
\rangle^2} + \frac{\langle G,H \rangle}{\langle H,H \rangle}\dot{H}.
\end{align}
There is a constant $a>0$ such that $\|G\|_\infty,\|H\|_\infty\le a$
for all $s$. Furthermore, $\|\dot{G}\|_\infty\le
c_{k+1,r}\|\dot{U}_s\|_\infty$ and $\|\dot{H}\|_\infty\le
c_{r+1,0}\|\dot{U}_s\|_\infty$ by induction. Since
\[
\Bigl|\frac{d}{ds}\langle H, H\rangle\Bigr| = 2 |\langle H,\dot{H}
\rangle| \le 2\|H\|_\infty\|\dot{H}\|_\infty < b \|\dot{U}_s\|_\infty
\]
for some constant $b>0$, it follows that $\langle H,H\rangle$
is a Lipschitz function of $U$ (as a
real-valued function), and hence if $\UU_{k+1,r+1}$ is small enough,
then $\langle H,H\rangle \ge c > 0$ for all $U\in\UU_{k+1,r+1}$.
Hence \eqref{DPSI} implies that $c_{k+1,r+1}$ exists, which proves
the Claim.

By Lemma \ref{LEM:NONFIN-DIFF}, we have $\langle \Phi(W), \Phi(W)
\rangle=\langle \Phi(W), t(F_{m+1}^\ddag,W) \rangle=0$, which is a
contradiction.
\end{proof}

In particular, it follows that the functions $t_2(F^\ddag,W)$ are
linearly dependent, which implies that finitely forcible graphons are
in a sense ``homogeneous''.

\begin{corollary}\label{COR:FIN-NEC}
Let $W\in\WW$ be finitely forcible. Then there is a nonzero simple
2-labeled connected quantum graph $f$ such that $t_2(f,W)=0$ almost
everywhere.
\end{corollary}

\begin{proof}
As remarked before, every finitely forcible graphon can be forced by
connected graphs $F_1,\dots,F_m$. Then either the functions
$t_2(F_1^\ddag,W),\dots, t_2(F_m^\ddag,W)$ are linearly dependent, or
else they span $\LL(W)$, and so for any simple connected graph
$F_{m+1}$ they span $t_2(F_{m+1}^\ddag,W)$. In either case we get a
finite set $F_1,\dots,F_k$ of (distinct) connected graphs such that
the functions $t_2(F_1^\ddag,W),\dots, t_2(F_k^\ddag,W)$ are linearly
dependent. Suppose that $\sum_i a_i t_2(F_i^\ddag,W) = 0$, then
$t_2(f,W)=0$ for $f=\sum_i F_i^\ddag$.

It is clear that $f$ is a simple connected 2-labeled quantum graph.
It is also clear that $f\not=0$: From each constituent of $F_i^\ddag$
we can reconstruct $F_i$ by connecting the labeled nodes by an edge
and deleting the labels. Hence constituents coming from different
$F_i^\ddag$ are different and they cannot cancel each other.
\end{proof}

The last corollary can be used to show that ``most'' graphons are not
finitely forcible.

\begin{theorem}\label{THM:DENSE}
The set of finitely forcible graphons is of first category in
$L_2([0,1]^2)$.
\end{theorem}

\begin{proof}
We claim that for a fixed set $\{F_1,\dots,F_k\}$ of connected simple
2-labeled graphs, the set of graphons $W$ for which there is a
nonzero quantum graph $f=\sum_{i=1}^k a_iF_i$ composed of these $F_i$
satisfying an equation $t_2(f,W)=0$ is nowhere dense. Let us fix a
$W$, we want to show that an arbitrary neighborhood of $W$ contains a
graphon $W'$ such that $t_2(F_1,W'),\dots,t_2(F_k,W')$ are linearly
independent. This will be enough, since $t_2$ is continuous and so
there is an open set $\UU$ in the neighborhood such that
$t_2(F_1,U),\dots,t_2(F_k,U)$ are linearly independent for all
$U\in\UU$.

Lemma 5 of \cite{ELS} implies that there are graphons $U_1,\dots,U_k$
such that the matrix $(t(F_i,U_j))_{i,j=1}^k$ is nonsingular. We may
assume that $\|W\|_\infty,\|U_1\|_\infty,\dots,\|U_k\|_\infty\le 1$.
For $0<\eps<1/k$, define $W^\eps=(1-k\eps)W\oplus(\eps) U_1\oplus
\dots\oplus(\eps) U_k$ (so the components of $W^\eps$ are
$W,U_1,\dots,U_k$, scaled by $1-k\eps,\eps,\dots,\eps$).

First we show that $W^\eps\to W$ in $L_2[0,1]^2$ if $\eps\to0$. Let
$W_\eps=(1-k\eps)W\oplus (k\eps) 0$. Then
\[
\|W^\eps - W_\eps\|_2^2 =
\eps^2(\|U_1\|_2^2+\cdots+\|U_k\|_2^2)\longrightarrow0 \qquad
(\eps\to 0),
\]
so it suffices to show that $\|W - W_\eps\|_2\to 0$. This is easy if
$W$ is a stepfunction with interval steps, and it follows for general
$W$ as these can be approximated by such stepfunctions in $L_2$.

Now suppose that $t_2(F_1,W^\eps),\dots,t_2(F_k,W^\eps)$ are linearly
dependent, so that there are real numbers $a_i$ such that
\[
\sum_{i=1}^k a_i t_2(F_i,W^\eps)(x,y)=0
\]
for all $x,y\in[0,1]$. If we integrate only over the points in in the
interval $(1-k\eps+(j-1)\eps,1-k\eps+j\eps)$, we get that
\[
\sum_{i=1}^k a_i \eps^{|V(F_i)|} t(F_i,U_j)=0 \qquad (j=1,\dots,k)
\]
(here we use that every connected component of each $F_i$ contains a
labeled node). But this contradicts the nonsingularity of the matrix
$(t(F_i,U_j))_{i,j=1}^k$.
\end{proof}

If $W$ is only finitely forced in $\WW_0$, then we get a weaker
condition:

\begin{lemma}\label{LEM:NONFIN-DIFF-2}
Suppose that $W\in\WW_0$ is forced in $\WW_0$ by the simple graphs
$F_1,\dots,F_m$, where we include $F_1=K_1$. Also suppose that there
exists an open ball $\UU\subseteq\WW$ about $W$ and a Lipschitz map
$\Phi:~\UU\to L_\infty[0,1]$ such that for all $U\in\UU$, $\langle
\Phi(U), t_1(F_i^\dag,U) \rangle=0$ $(i=1,\dots,m)$. Then $\langle
\Phi(W), t_1(F^\dag,W) \rangle=0$ for every simple graph $F$.
\end{lemma}

\begin{proof}
Since we are not using this lemma, we only sketch the proof: instead
of changing the values of the function $W$, we add a weight function.
Let $\alpha\in L_\infty[0,1]$ be a function with
$\int\limits_0^1\alpha(y)\,dy=1$. For every multigraph $F$, we define
\[
t(F,W,\alpha)=\int\limits_{[0,1]^{V(F)}} \prod_{i\in V(F)}
\alpha(y_i) \prod_{ij\in E(F)}W(y_i,y_j)\prod_{i\in V(F)} dy_i,
\]
and for every $1$-labeled graph $F$, we define
\[
t_1(F,W,\alpha)(y_1)=\int\limits_{[0,1]^{V(F)\setminus \{1\}}}
\prod_{i\in V(F)} \alpha(y_i) \prod_{ij\in E(F)}W(y_i,y_j)\prod_{i\in
V(F)\setminus \{1\}} dy_i.
\]
Through the change of variables $x=\int\limits_0^y\alpha(z)\,dz$, it
is easy to construct a graphon $W_\alpha\in \WW_0$ for which
$t(F,W,\alpha)=t_F(F,W_\alpha)$ for every multigraph $F$. Note that
in particular $t(K_1,W,\alpha)=\int\limits_0^1\alpha(y)\,dy=1$.
Furthermore,
\[
t_1(F,W_\alpha)\bigl(\int\limits_0^y\alpha(z)\,dz\bigr)=
t_1(F,W,\alpha)(y)
\]
for every $1$-labeled graph $F$.

Let $(\alpha_s:~0\le s\le 1)$ be a differentiable family of functions
in $L_\infty[0,1]$ such that $\alpha_0\equiv1$ and
$\int\limits_0^1\alpha_s(y)\,dy=1$. Lemma \ref{LEM:WDIFF} can be
replaced by the formula
\[
\frac{d}{ds} t(F,W,\alpha_s) = \Bigl\langle \frac{\partial}{\partial
s}\alpha_s, t_1(F^\dag,W,\alpha_s) \Bigr\rangle.
\]
The proof is concluded by solving the differential equation
\[
\frac{\partial}{\partial s} \alpha_s(y) =
\Phi(W_{\alpha_s})\bigl(\int\limits_0^y \alpha_s(z)\, dz\bigr),
\]
in the Banach space $L_\infty[0,1]$, similarly as in the proof of
Lemma \ref{LEM:NONFIN-DIFF}.
\end{proof}

\begin{corollary}\label{PROP:FIN-NEC-2}
Let $W\in\WW_0$ be finitely forcible in $\WW_0$. Then there is a
simple connected nonzero 1-labeled quantum graph $f\not=0$ such that
$t_1(f,W)=0$ almost everywhere.
\end{corollary}

\section{Open Problems and further directions}

It does not seem easy to characterize finitely forcible functions.
Let us offer a few conjectures. The next question might be easy but
the examples and theorems in the present paper don't answer it.

\begin{question}
Is there a non-constant continuous (or smooth) fnction on $[0,1]^2$
which is finitely forcible? (As we have seen, the simplest
candidates, namely polynomial functions, don't work.)
\end{question}

We believe that in Theorem \ref{MON-POLY}, the assumption that $p$ is
monotone can be omitted:

\begin{conjecture}\label{SUFF}
For every symmetric 2-variable polynomial $p$, the function
$\one_{p(x,y)\ge 0}$ is finitely forcible in $\WW$. (Using ad hoc
tricks, the proof given in Section \ref{POLYN} can be extended to
some non-monotone polynomials, for example, to $(1/2-x-y)(3/2-x-y)$.)
\end{conjecture}

We can try to generalize the results of Section \ref{POLYN} to more
variables. Here is an interesting special case:

\begin{question}\label{QU:SPHERE}
Is the following graphon finitely forcible: the underlying
probability space is the uniform distribution on the surface of the
unit sphere $S^2$, and $W(x,y)=1$ if $x$ and $y$ are closer than
$90^\circ$, and $W(x,y)=0$ otherwise?
\end{question}

It is not clear whether the two notions of forcibility we have
considered are really different.

\begin{question}
If a function $W\in\WW_0$ is finitely forcible in $\WW_0$, is it also
finitely forcible in $\WW$?
\end{question}

We don't know too much about algebraic operations which generate new
forcible functions. For example, it is unreasonable to expect that
the sum of two forcible functions is forcible, since the sum depends
on the concrete representation of the graphons (not just on their
weak isomorphism types). However the next question is natural.

\begin{question}\label{QU:TENSOR}
Is the tensor product $U\otimes W$ of two finitely forcible graphons
$U$ and $W$ forcible?
\end{question}

Corollary \ref{UB-FORCE} suggests the following problem:

\begin{question}
For which finite graphs $G$ is the infinite lexicographic power of
$G$ finitely forcible?
\end{question}

Our motivation for the study of finitely forcible graphons was to
understand the structure of extremal graphs. This would be fully
justified by the following conjecture:

\begin{conjecture}\label{CONJ:UNIQUE}
If a finite set of constraints of the form $t(F_i,W)=a_i$
($i=1,\dots,k$) is satisfied by some graphon, then it is satisfied by
a finitely forcible graphon. This conjecture would imply the
(imprecise) fact that {\it every extremal graph problem has a
finitely forcible solution}.
\end{conjecture}

Let us state the problem mentioned in Remark \ref{INFINITESIMAL}:

\begin{question}\label{CONJ:INF-FF}
Is every infinitesimally finitely forcible graphon also finitely
forcible (and hence, a stepfunction)?
\end{question}

The topology of the set $T(W)$ introduced in Section \ref{TYPICAL}
gives rise to some interesting problems. It is easy to see that
$R(W)\cap T(W)$ is dense in $T(W)$ (in the topology of $L_1[0,1]$,
and if two graphons $W$ and $U$ are weakly isomorphic then $T(U)$ is
homeomorphic to $T(W)$.

Surprisingly, in each of the finitely forcible examples of this paper
$T(W)$ is a finite dimensional compact topological space. For
positive supports of monotone polynomials, $T(W)$ is homeomorphic
with the interval $[0,1]$. The topology of the regular CR-graphon
corresponding to the binary tree is the Cantor set $\{0,1\}^\infty$.
The examples constructed in Section \ref{SEC:IRR} correspond to the
one-point compactification of the natural numbers.

This topological space was introduced and studied in \cite{LSz2,LSz7},
where it was proved that if $t(F,W)=0$ for some signed bipartite
graph $F$, then $T(W)$ is finite dimensional and compact. The
following conjectures would lead to the same conclusion from a
different assumption.

\begin{conjecture}
If $W$ is finitely forcible in $\WW_0$ then $T(W)$ is a compact
space. (We can't even prove that $T(W)$ is locally compact.)
\end{conjecture}

\begin{conjecture}\label{CONJ:DIM}
If $W$ is finitely forcible then $T(W)$ is finite dimensional. (We
intentionally do not specify which notion of dimension is meant
here---a result concerning any variant would be interesting.) Note
that Corollary \ref{INF-RANK} implies that the linear hull of $T(W)$
is infinite dimensional unless $T(W)$ is a finite set.

In our examples $T(W)$ is either $0$-dimensional of $1$-dimensional.
This is probably due to the fact that we have only found very simple
examples.
\end{conjecture}

\begin{question}\label{highdim}
Is there a finitely forcible graphon $W$ such that $T(W)$ is
homeomorphic with $[0,1]^2$? (A positive answer would follow from a
positive answer to question \ref{QU:TENSOR}.)
\end{question}

One can also consider a more direct notion of dimension. We define
the {\it dimension of the graphon} $W$ as the infimum of all $c>0$
such that for every $\eps>0$ there is a stepfunction $W_\eps$ with
$O((1/\eps)^c)$ steps such that $\|W-W_\eps\|_\square\leq\eps$. It
was shown in \cite{LSz2}) that the dimensions $W$ and $T(W\circ W)$
are related.

The dimension of $W$ can be described in terms of the number of
classes in weak Szemer\'edi partitions (introduced by Freeze and
Kannan \cite{FK}). So a positive answer to Conjectures
\ref{CONJ:UNIQUE} and \ref{CONJ:DIM} would imply that extremal graph
problems have solution with efficient (polynomial-size) weak
Szemer\'edi partitions. This could explain (in a weak sense) why
Szemer\'edi partitions are so important in extremal graph theory.

\subsection*{Acknowledgement}

Our thanks are due to Vera S\'os and to Hamed Hatami for their
interest in this work. We also thank them and the anonymous referees
for their questions and many suggestions that lead to substantial
improvements.

\end{document}